\documentclass[11pt, leqno, a4paper]{article}

\usepackage[top=2.6cm, bottom=2.6cm, left=2.5cm, right=2.5cm]{geometry}

\usepackage{latexsym,amssymb}
\usepackage{amsmath,amsthm}
\usepackage{color} 
\usepackage{mathtools} 

\usepackage{multirow,eepic} 

\newtheorem{Theorem}{\bf Theorem}[section]
\newtheorem{Lemma}{\bf Lemma}[section]
\newtheorem{Proposition}{\bf Proposition}[section]

\newtheorem{Remark}{\bf Remark}[section]

\allowdisplaybreaks[4]

\numberwithin{equation}{section}

\begin{document}
\title{Existence and nonexistence of solutions for \\ the heat equation with a superlinear source term} 
\author{
        Yohei Fujishima\footnote{e-mail address: fujishima@shizuoka.ac.jp} \\ \\ 
        {\small Department of Mathematical and Systems Engineering} \\ 
        {\small Faculty of Engineering, Shizuoka University} \\ 
        {\small 3-5-1 Johoku, Hamamatsu 432-8561, Japan} \\ \\ 
        Norisuke Ioku\footnote{e-mail address: ioku@ehime-u.ac.jp} \\ \\ 
        {\small Graduate School of Science and Engineering, Ehime University} \\ 
        {\small Matsuyama, Ehime 790-8577, Japan} 
        } 
\date{} 
\pagestyle{myheadings}
\markboth{Y. Fujishima and N. Ioku}{Existence results for nonlinear heat equations}

\maketitle 
\begin{abstract} 
Classification theory on the existence and non-existence of local in time solutions 
for initial value problems of nonlinear heat equations are investigated. 
Without assuming a concrete growth rate on a nonlinear term, 
we reveal the threshold integrability of initial data which classify existence and nonexistence of solutions via a quasi-scaling and its invariant integral.
Typical nonlinear terms, for instance polynomial type, exponential type and its sum, product and composition, can be treated as applications. 
%
\end{abstract} 
\vspace{5pt} 


\section{Introduction} 
\label{section:1} 
We consider existence and nonexistence of solutions for a heat equation with general nonlinearity 
\begin{equation} 
\label{eq:1.1} 
\left\{ 
\begin{array}{lll} 
\partial_t u=\Delta u+f(u) &\mbox{in} &{\mathbb{R}}^N\times (0,T), \vspace{3pt} \\ 
u(x,0)=u_0(x)\ge 0 &\mbox{in} &{\mathbb{R}}^N, 
\end{array} 
\right. 
\end{equation} 
where $\partial_t=\partial/\partial t$, $N\ge 1$, $T>0$, $u_0$ is a nonnegative measurable initial function and 
$f\in C^1([0,\infty))$ is a positive monotonically increasing function in $(0,\infty)$, that is, 
\begin{equation} 
\label{eq:1.2} 
f(s)>0, \quad f'(s)>0 \quad\mbox{for all}\,\,\, s\in (0,\infty). 
\end{equation} 
In the following, for suitable Banach space $X$, 
we say a function $u\in C^{2,1}(\mathbb{R}^N\times (0,T))$ is a classical solution in $X$ for problem~\eqref{eq:1.1} if 
$u$ satisfies the equation in the classical sense
and
$\|u(t)-e^{t\Delta}u_0\|_{X}\to 0$ as $t\to 0$, 
where $e^{t\Delta}u_0$ denotes the solution of the heat equation with the initial data $u_0$. 

It follows from the standard argument that problem~\eqref{eq:1.1} possesses the unique classical solution in $L^\infty(\mathbb{R}^N)$ 
for general nonlinearities $f\in C^1([0,\infty))$ if $u_0\in L^\infty(\mathbb{R}^N)$.  
On the other hand, for the case $u_0\not\in L^\infty({\mathbb R}^N)$, 
existence results of solutions for problem~\eqref{eq:1.1} heavily depend on the growth rate of the nonlinear term 
$f$. 
One of the typical examples of $f$ is a power type nonlinearity, that is, 
\begin{equation} 
\label{eq:1.3} 
\partial_t u=\Delta u+u^p, \quad x\in {\mathbb{R}}^N, \,\,\, t>0, \qquad  
u(x,0)=u_0(x)\ge 0, \quad x\in \mathbb{R}^N, 
\end{equation} 
where $p>1$. 
This equation for unbounded initial data has been studied intensively since the pioneering work due to Weissler \cite{W1}, 
and there hold the following: 
\begin{itemize} 
\item 
if $r\ge \frac{N}{2}(p-1)$ and $r>1$ or $r>\frac{N}{2}(p-1)$ and $r\ge 1$,
then for any $u_0\in L^r(\mathbb{R}^N)$, there exist a constant $T>0$ and 
a local in time classical solution $u\in C([0,T];L^r(\mathbb{R}^N))$ for problem~\eqref{eq:1.3}. 
\item 
if $\frac{N}{2}(p-1)>1$ and $1\le r<\frac{N}{2}(p-1)$,
then there exists an initial function $u_0\in L^r(\mathbb{R}^N)$ such that problem~\eqref{eq:1.3} can not possess 
any local in time nonnegative classical solutions. 
\end{itemize}  
See \cite{AD}, \cite{BP}, \cite{BC}, \cite{G}--\cite{IKS}, \cite{LRSV}, \cite{NS}, \cite{RS}--\cite{W1} for 
existence and nonexistence of solutions for nonlinear parabolic equations and qualitative properties of solutions. 
We also refer to \cite{LSU} and \cite{QS}, which include good references concerning parabolic equations. 
It is remarkable that the critical exponent $r_c:=\frac{N}{2}(p-1)$, 
which gives the classification of existence and nonexistence of solutions for problem~\eqref{eq:1.3},
arises from the scale invariant property of~\eqref{eq:1.3}
under the scaling transformation
%
\begin{equation} 
\label{eq:1.4} 
u_\lambda(x,t):=\lambda^\frac{2}{p-1}u(\lambda x,\lambda^2 t), \qquad \lambda>0,
\end{equation}
in the sense that
$\|u_\lambda(\cdot, 0)\|_{L^r(\mathbb{R}^N)}=\|u(\cdot, 0)\|_{L^r(\mathbb{R}^N)}$ if and only if $r=r_c$. 
The above assertions due to Weissler show that the critical exponent $r=r_c$
can be a priori found by the scaling transformation and its invariant norm of problem~\eqref{eq:1.3}.
However, in the case of general nonlinearity $f(u)$, 
it is not clear which integral should be controlled for the classification of existence and nonexistence of solutions,
since problem~\eqref{eq:1.1} 
does not possess a scale invariant property for general nonlinearity $f$.

%
%
%
In this paper, we reveal the threshold integrability of $u_0$ to classify existence and nonexistence of solutions for problem~\eqref{eq:1.1} 
without any concrete assumptions on the behavior of $f=f(s)$ near $s=\infty$.
To this end, 
we apply the ``quasi'' scaling proposed by the first author of this paper in \cite{F}:
\begin{equation} 
\label{eq:1.6} 
u_\lambda(x,t):=F^{-1}\Big[ 
\lambda^{-2} F(u(\lambda x, \lambda^2 t)) 
\Big], \qquad \lambda>0, 
\end{equation} 
where 
\begin{equation}\notag
F(s):=\int_s^\infty \frac{du}{f(u)}
\end{equation} 
and 
$F^{-1}$ is the inverse function of $F$.
We mention that
the transformation \eqref{eq:1.6} does not preserve the equation \eqref{eq:1.1},  
however, the main term of \eqref{eq:1.1} are unchanged. 
In fact, 
for the solution $u$ of \eqref{eq:1.1}, 
the 
function $u_{\lambda}$ defined by \eqref{eq:1.6}
satisfies
\begin{equation} 
\label{eq:1.7} 
\partial_t u_\lambda=\Delta u_\lambda+f(u_\lambda)+f(u_\lambda)f(u)^{-2}F(u)^{-1}|\nabla u|^2\Big[ 
f'(u) F(u)-f'(u_\lambda)F(u_\lambda) 
\Big].
\end{equation} 
It should be emphasized that this striking 
transformation \eqref{eq:1.6} is a generalization of \eqref{eq:1.4} for polynomial nonlinearity
since \eqref{eq:1.6} coincides with \eqref{eq:1.4} if $f(u)=u^p$ and the remainder term in \eqref{eq:1.7} becomes zero.
%

We now introduce 
a significant property which arises from the transformation~\eqref{eq:1.6}, 
that is, the scale invariant property 
\begin{equation} 
\label{eq:1.8} 
\int_{\mathbb{R}^N}F(u_\lambda(x,0))^{-\frac{N}{2}}\, dx=\int_{\mathbb{R}^N}F(u(x,0))^{-\frac{N}{2}}\, dx, 
\qquad \lambda>0. 
\end{equation} 
This property plays an essential role for classifying existence and nonexistence as $L^{r_c}$ norm for problem \eqref{eq:1.4}.
Furthermore, we focus on the  limit of $f'(s)F(s)$ as $s\to \infty$:
$$
A:=\lim_{s\to\infty}f'(s)F(s),
$$
since 
the behavior of the function $f'F$ controls
the remainder term of \eqref{eq:1.7}.
%
%
Throughout this paper, we assume that the above limit $A$ always exists 
and 
$f$ is 
superlinear in the sense that 
\begin{equation} 
\label{eq:1.5} 
F(s)=\int_s^\infty \frac{du}{f(u)}<\infty 
\end{equation} 
for all $s>0$. 
Note that the limit $A$ always exist for typical examples of $f$ such as 
\begin{equation} 
\notag 
f(u)=u^p \ (p>1), \quad u^p+u^q \ (p>q>1), \quad e^u, \quad e^{u^2}. 
\end{equation} 

The main purpose of this paper is to show that the integrability driven from \eqref{eq:1.8} 
implies the classification theorems of existence of solutions for problem~\eqref{eq:1.1} with general nonlinearity $f$. 
%
To state the results, we introduce some notation. 
For $x\in \mathbb{R}^N$ and $\rho>0$, we denote by $B_\rho(x)$ the ball of radius $\rho$ centered at $x$.  
For $1\le p<\infty$, define the uniformly local $L^p$ space $L^p_{ul,\rho}(\mathbb{R}^N)$ 
by 
$$
 L^p_{ul,\rho}(\mathbb{R}^N)
 := \left\{
  u\in L^p_{loc}(\mathbb{R}^N):\|u\|_{L^p_{ul,\rho}(\mathbb{R}^N)}:=\sup_{y \in \mathbb{R}^N} 
  \left(\int_{B_{\rho}(y)}|u(x)|^p\, dx\right)^{1/p}<\infty
 \right\}. 
$$ 
We denote by $\mathcal{L}^p_{ul,\rho}(\mathbb{R}^N)$ the closure of the space of bounded uniformly continuous functions 
$BUC(\mathbb{R}^N)$ in the space $L^p_{ul,\rho}(\mathbb{R}^N)$, that is, 
$$
\mathcal{L}^p_{ul,\rho}(\mathbb{R}^N) 
:=\overline{BUC(\mathbb{R}^N)}^{\|\cdot\|_{L^p_{ul,\rho}(\mathbb{R}^N)}}. 
$$
Define 
$$
(e^{t\Delta}\varphi)(x):=(4\pi t)^{-\frac{N}{2}}\int_{\mathbb{R}^N}e^{-\frac{|x-y|^2}{4t}}\varphi(y)\, dy 
$$
for $\varphi\in L_{ul,\rho}^1(\mathbb{R}^N)$. 
Then $e^{t\Delta}\varphi$ gives the solution for the heat equation with the initial data $\varphi$. 
\vspace{5pt} 

We are ready to state our main results. 
We first state the existence results of local in time solutions for problem~\eqref{eq:1.1}. 
\begin{Theorem}
\label{Theorem:1.1} 
Let $N\ge 1$, $\rho>0$ and $f\in C^1([0,\infty))$ satisfy \eqref{eq:1.2} and \eqref{eq:1.5}. 
Assume that the limit $A=\displaystyle\lim_{s\to\infty}f'(s)F(s)$ exists and that there exists a constant $s_1>0$ such that  
\begin{equation} 
\label{eq:1.9} 
f'(s)F(s)\le A \quad\mbox{for all}\,\,\, s\ge s_1. 
\end{equation} 
{\rm (i) (Subcritical case)} 
Let $r>0$ be a constant satisfying 
\begin{equation} 
\label{eq:1.10} 
r\ge A-1, \qquad r>\frac{N}{2}, 
\end{equation} 
and assume that a nonnegative initial function $u_0$ is measurable and satisfies 
\begin{equation} 
\label{eq:1.11} 
F(u_0)^{-r}\in L_{ul,\rho}^1(\mathbb{R}^N). 
\end{equation} 
\begin{itemize}
\item 
If $A>1$,
there exist $T>0$ and a local in time classical solution $u\in C^{2,1}(\mathbb{R}^N\times (0,T))$ of \eqref{eq:1.1} satisfying 
the initial value problem in the following sense: 
\begin{equation}
 \label{eq:1.12}
 \lim_{t\to 0}\|u(t)-e^{t\Delta}u_0\|_{L^{\frac{r}{A-1}}_{ul,\rho}(\mathbb{R}^N)}
 =0.
 \end{equation}
Furthermore, the existence time $T$ can be estimated to satisfy 
\begin{equation}
\label{eq:1.13} 
\begin{split} 
& T^{\frac{N}{2}(1-\frac{1}{A})}\rho^{-N(1-\frac{1}{A})} \\ 
& \quad 
+\max\left\{\|F(u_0)^{-r}\|_{L^1_{ul,\rho}(\mathbb{R}^N)}, F(s_1)^{-r}\rho^N\right\} 
\left(T^{r-\frac{N}{2A}}\rho^{-N\cdot\frac{A-1}{A}}+T^{r-\frac{N}{2}}\right) 
\ge 
\gamma  
\end{split} 
\end{equation}
where $\gamma>0$ depends only on $N$, $A$ and $r$. 
\item 
If $A=1$, 
there exist $T>0$ and a local in time classical solution $u\in C^{2,1}(\mathbb{R}^N\times (0,T))$ of \eqref{eq:1.1} satisfying 
the initial value problem in the following sense: 
\begin{equation}
 \label{eq:1.14}
 \lim_{t\to 0}\|u(t)-e^{t\Delta}u_0\|_{L^{\infty}(\mathbb{R}^N)}
 =0.
\end{equation}
Furthermore, the existence time $T$ can be estimated to satisfy 
\begin{equation}
 \label{eq:1.15} 
 \begin{split} 
 & T^{\frac{\epsilon N}{2}}\rho^{-\epsilon N} \\ 
 & + \gamma_\epsilon \max\left\{\|F(u_0)^{-r}\|_{L^1_{ul,\rho}(\mathbb{R}^N)}, F(s_1)^{-r}\rho^N\right\} 
 \left(T^{r-\frac{N}{2}(1-\epsilon)}\rho^{-\epsilon N}+T^{r-\frac{N}{2}}\right) 
 \ge 
 \gamma 
 \end{split} 
\end{equation}
for any sufficiently small $\epsilon>0$, 
where $\gamma$ depends only on $N$ and $r$, 
and $\gamma_\epsilon$ is a positive constant depending only on $N$, $r$ and $\epsilon$ 
satisfying
$\gamma_\epsilon\to \infty $ as $\epsilon\to 0$. 
\end{itemize}
\noindent
{\rm (ii) (Critical case)} 
Let 
\begin{equation} 
\label{eq:1.16} 
r=\frac{N}{2}>A-1, 
\end{equation} 
and assume that a nonnegative initial function $u_0$ is measurable and satisfies 
\begin{equation} 
\label{eq:1.17} 
F(u_0)^{-r}\in {\cal L}_{ul,\rho}^{1}({\mathbb R}^N). 
\end{equation} 
Then there exist $T>0$ and a local in time classical solution $u\in C^{2,1}(\mathbb{R}^N\times (0,T))$ 
of \eqref{eq:1.1} satisfying \eqref{eq:1.12} for the case $A>1$ and \eqref{eq:1.14} for the case $A=1$. 
\end{Theorem}
In view of Theorem~\ref{Theorem:1.1}, 
the balance of two important factors, 
the behavior of the function $f=f(s)$ as $s\to\infty$ and the singularity of $u_0$ (given in \eqref{eq:1.11} or \eqref{eq:1.17}),
is controlled by the constant $\displaystyle A=\lim_{s\to \infty}f'(s)F(s)$ via the function $f$. 
Note that it relies on the setting of $f$ and the space of the initial data whether the solution $u(t)$ converges to $u_0$ as $t\to 0$, 
so we do not discuss this problem under general setting of Theorem~\ref{Theorem:1.1}. 
However, we consider the convergence of the solution to the initial data for concrete examples of $f$ in Section~\ref{section:5}. 
See also Remark~\ref{Remark:1.3}. 

\begin{Remark} 
\label{Remark:1.1}
It must hold that $A\ge 1$ for all superlinear functions $f\in C^1([0,\infty))$ satisfying \eqref{eq:1.2} and \eqref{eq:1.5}. 
Assume that $A<1$. 
Then there exist $s_0>0$ and $\epsilon\in (0,1)$ such that
$f'(s)F(s)\le 1-\epsilon$ for all $s\ge s_0$. 
Since $f'(s)=F''(s)/F'(s)^2$ and $F'(s)=-1/f(s)<0$, it holds 
$$
 \frac{F''(s)}{F'(s)}\ge (1-\epsilon) \frac{F'(s)}{F(s)}
$$
for all $s\ge s_0$. 
Integrating both sides of above inequality on $(s_0,s)$, we have
$$
-F'(s)\ge CF(s)^{1-\epsilon}, 
$$
where $C>0$ is a constant.
This together with a simple calculation 
contradicts the positivity of $F$. 
\end{Remark}
\begin{Remark} 
\label{Remark:1.2} 
One can easily check that 
$$
\lim_{s\to\infty}f'(s)F(s)=\left\{ 
\begin{array}{cll} 
\dfrac{p}{p-1} &\mbox{if} &f(u)=u^p\ \mbox{with}\ p>1, \vspace{3pt} \\ 
1 &\mbox{if} &f(u)=e^u\ \mbox{or}\ f(u)=e^{u^2}. 
\end{array} 
\right. 
$$
Therefore the constant $A$ becomes small if the growth rate of $f$ is rapid. 
\end{Remark} 
\begin{Remark} 
\label{Remark:1.3} 
Under the condition~\eqref{eq:1.9}, 
as in Remark~{\rm\ref{Remark:1.1}}, we have 
$$
\frac{F''(s)}{F'(s)}\ge A\frac{F'(s)}{F(s)}, 
$$
which implies that $-F'(s)\ge CF(s)^A$ for all sufficiently large $s>0$. 
Here $C>0$ is a constant independent of $s$. 
Then, by a simple calculation we can check that 
$s
\le CF(s)^{-(A-1)}$ for some $C>0$ and all sufficiently large $s>0$ and that 
$F(u_0)^{-r}\in L_{ul, \rho}^1(\mathbb{R}^N)$ yields 
$u_0\in L_{ul,\rho}^{\frac{r}{A-1}}(\mathbb{R}^N)$ for the case $A>1$ 
provided that either \eqref{eq:1.10} and \eqref{eq:1.11} or \eqref{eq:1.16} and \eqref{eq:1.17} hold. 
Therefore, the convergence \eqref{eq:1.12} of $u(t)$ 
is reasonable from the viewpoint of the singularity of $u_0$. 
On the other hand, \eqref{eq:1.11} does not imply that 
$u_0\in {\cal L}_{ul,\rho}^{\frac{N}{2}\cdot\frac{1}{A-1}}(\mathbb{R}^N)$ for the case $A>1$ in general. 
So we do not know whether $e^{t\Delta}u_0$ can be replaced by $u_0$ in \eqref{eq:1.12}. 
See also Lemma~{\rm \ref{Lemma:2.2}}. 
\end{Remark} 

We next state the nonexistence results of local in time solutions for problem~\eqref{eq:1.1}.
%
\begin{Theorem}[\bf Nonexistence] 
\label{Theorem:1.3}
Let $N\ge 1$, $\rho>0$ and $f\in C^1([0,\infty))$ satisfy \eqref{eq:1.2} and \eqref{eq:1.5}. 
Assume that the limit $A=\displaystyle\lim_{s\to\infty}f'(s)F(s)$ exists. 
Furthermore, assume that $A-1<N/2$ and there exists a constant $s_2>0$ such that  
\begin{equation} 
\label{eq:1.20} 
f'(s)F(s)\ge A \quad\mbox{for all}\,\,\, s\ge s_2. 
\end{equation} 
Then, for any $r\in [A-1, N/2)$ if $A>1$ or any $r\in (0,N/2)$ if $A=1$, there exists a nonnegative measurable initial function $u_0$ satisfying 
$F(u_0)^{-r} \in {L}_{ul,\rho}^1(\mathbb{R}^N)$ 
such that there can not exist nonnegative classical solutions of \eqref{eq:1.1} satisfying 
the initial value problem in the sense \eqref{eq:1.12} or \eqref{eq:1.14}. 
\end{Theorem} 

In some examples of $f$, for instance $f(u)=u^p+u^q\ (p>q>1)$,
the condition~\eqref{eq:1.20} does not hold. 
However, for the case $A>1$, it seems possible 
to avoid this difficulty by considering some approximation of $f$ from above (See Theorem~\ref{Theorem:5.1} and its proof). 
On the other hand, for the case $A=1$, 
similar calculations as in Remark~\ref{Remark:1.3} can be carried out, 
and the condition~\eqref{eq:1.20} with $A=1$ implies that 
there exists a constant $C>0$ such that $f(u)\le e^{Cu}$. 
Therefore, Theorem~\ref{Theorem:1.3} is no longer available for rapidly growing nonlinearity such as $f(u)=e^{u^2}$. 
However, even for the case $f$ has rapid growth, that is, $A=1$ 
and $f'(s)F(s)$ converges to $1$ from below, 
there holds the following nonexistence result.
%
\begin{Theorem}[\bf Nonexistence]  
\label{Theorem:1.4} 
Let $N\ge 1$, $\rho>0$ and $f\in C^1([0,\infty))$ satisfy \eqref{eq:1.2} and \eqref{eq:1.5}. 
Assume that $A=1$ and 
there exists a constant $s_3>0$ such that 
$$	 
f'(s)F(s)\le 1 \quad\mbox{for all}\,\,\, s\ge s_3. 
$$
Then, for any $r\in (0, N/2)$, 
there exists a nonnegative measurable initial function $u_0$ satisfying 
$F(u_0)^{-r}\in L^1_{ul, \rho}(\mathbb{R}^N)$ such that there can not exist nonnegative classical solutions 
of \eqref{eq:1.1} satisfying the initial value problem in the sense \eqref{eq:1.14}. 
\end{Theorem}

So far as the authors know, 
these are the first results to characterize existence and nonexistence of solutions for problem~\eqref{eq:1.1} 
without any assumptions on the growth rate of the nonlinear term $f$. 
\vspace{5pt} 

As an application of our main results, we treat the case $f(u)=e^{u^2}$. Several other applications are considered in Section~\ref{section:5}.
Consider 
\begin{equation} 
\label{eq:1.21} 
\left\{ 
\begin{array}{ll} 
\partial_t u=\Delta u+e^{u^2}, &x\in \mathbb{R}^N, \,\,\, t>0, \vspace{3pt} \\ 
u(x,0)=u_0(x)\ge 0, &x\in \mathbb{R}^N. 
\end{array} 
\right. 
\end{equation} 
This type nonlinearity appears in view of the Trudinger-Moser inequality 
or the problem related to uniqueness results for \eqref{eq:1.3} with $N=2$,
and has been treated in \cite{IJMS}, \cite{I}, \cite{IRT} and \cite{RT}. 
See also \cite{NO}. 
Among others, in \cite{I}, \cite{IRT} and \cite{RT}, the authors discuss local in time existence 
and nonexistence of solutions for \eqref{eq:1.21} (with a slight modification on the nonlinearity)
in the Orlicz space $\operatorname{exp}L^2(\mathbb{R}^N)$. 
Here the Orlicz space $\operatorname{exp}L^2(\mathbb{R}^N)$ is the set of all functions $u_0\in L^1_{loc}(\mathbb{R}^N)$ satisfying 
\begin{equation} 
\label{eq:1.22} 
\int_{\mathbb{R}^N}\left[ 
e^{\lambda |u_0(x)|^2}-1 \right] dx<\infty 
\end{equation} 
for some $\lambda>0$. 
Local in time existence of solutions for problem~\eqref{eq:1.21} has been shown for initial data $u_0$ 
satisfying \eqref{eq:1.22} with sufficiently large $\lambda>0$, 
and they also prove that \eqref{eq:1.21} can not possess local in time solutions for 
some $u_0$ satisfying \eqref{eq:1.22} with sufficiently small $\lambda>0$. 
It seems that the critical integrability for $u_0$ which classifies local in time existence and nonexistence of solutions for \eqref{eq:1.21} is not known yet. 
In the following theorem, we apply Theorems~\ref{Theorem:1.1} and \ref{Theorem:1.4} to obtain the critical integrability of initial data $u_0$.  
\begin{Theorem}   
\label{Theorem:1.5} 
Let $N\ge 1$ and $\rho>0$. 
\begin{itemize} 
\item[\rm (i)] 
{\rm (Subcritical case)} 
Let $r>N/2$. 
For any nonnegative measurable function $u_0$ satisfying $|u_0|^r e^{r|u_0|^2}\in {L}_{ul,\rho}^1(\mathbb{R}^N)$, 
there exists a local in time classical solution for problem~\eqref{eq:1.21} satisfying \eqref{eq:1.14}. 
Furthermore, if $|u_0|^r e^{r|u_0|^2}\in {\cal L}_{ul,\rho}^1(\mathbb{R}^N)$, then 
\begin{equation} 
\label{eq:1.23} 
\lim_{t\to 0}\, \sup_{y\in\mathbb{R}^N}\int_{B_{\rho}(y)}|u(x,t)-u_0(x)|^r e^{r|u(x,t)-u_0(x)|^2}\, dx=0. 
\end{equation} 
The existence time $T$ can be chosen to satisfy
$$
T^{\frac{\epsilon N}{2}}\rho^{-\epsilon N} 
+ 
\gamma_\epsilon \max\left\{\|u_0^r e^{ru_0^2}\|_{L_{ul,\rho}^1(\mathbb{R}^N)}, \rho^N\right\} 
\left(T^{r-\frac{N}{2}(1-\epsilon)}\rho^{-\epsilon N}+T^{r-\frac{N}{2}}\right) 
\ge 
\gamma 
$$
for any sufficiently small $\epsilon>0$, 
where $\gamma$ depends only on $N$ and $r$, 
and $\gamma_\epsilon$ is a positive constant depending only on $N$, $r$ and $\epsilon$ 
and satisfies $\gamma_\epsilon\to \infty $ as $\epsilon\to 0$. 
\item[\rm (ii)] 
{\rm (Critical case)} 
Let $u_0$ be a nonnegative measurable function satisfying $|u_0|^\frac{N}{2} e^{\frac{N}{2}|u_0|^2}\in {\cal L}_{ul,\rho}^1(\mathbb{R}^N)$. 
Then there exists a local in time classical solution for problem~\eqref{eq:1.21} satisfying \eqref{eq:1.14} 
and \eqref{eq:1.23} with $r=N/2$. 
\item[\rm (iii)] 
{\rm (Nonexistence)} 
Let $0<r<N/2$. 
Then there exists a nonnegative measurable data $u_0$ such that 
$|u_0|^r e^{r|u_0|^2}\in {L}_{ul,\rho}^1(\mathbb{R}^N)$ 
and problem~\eqref{eq:1.21} can not possess any local in time 
classical solution $u$ satisfying \eqref{eq:1.14}. 
\end{itemize} 
\end{Theorem} 

We sketch the outline of the proof of our main theorems. 
In order to argue existence of solutions for \eqref{eq:1.1} with general nonlinearity $f$, 
we introduce a generalization of the Cole-Hopf transformation. 
In the case $A>1$, let $u$ satisfy $\partial_t u=\Delta u+f(u)$ and put 
\begin{equation} 
\label{eq:1.24} 
v(x,t):=F(u(x,t))^{-(A-1)}. 
\end{equation} 
Then $v$ satisfies 
$$
\partial_t v=\Delta v+(A-1)v^{1+\frac{1}{A-1}}+(A-1)F(u)^{-A-1}f(u)^{-2}|\nabla u|^2\Big[ 
f'(u)F(u)-A \Big].  
$$
Using the transformation~\eqref{eq:1.24}, we construct a supersolution of problem~\eqref{eq:1.1}, 
then we obtain a solution for problem~\eqref{eq:1.1} by a monotone method. 
See Proposition~\ref{Proposition:2.1}. 
For construction of a supersolution of \eqref{eq:1.1}, 
we utilize the solution for a heat equation with power type nonlinearity, that is, 
\begin{equation} 
\label{eq:1.25} 
\partial_t v=\Delta v+(A-1)v^{\frac{A}{A-1}}, \,\,\, x\in \mathbb{R}^N,\,\,\, t>0, \quad  
v(x,0)=F(u_0(x))^{-(A-1)}, \,\,\, x\in \mathbb{R}^N. 
\end{equation} 
Define
$$
\overline{u}(x,t):=F^{-1}\left(v(x,t)^{-1/(A-1)}\right). 
$$
Then $\overline{u}$ satisfies $\partial_t\overline{u}\ge \Delta\overline{u}+f(\overline{u})$, 
provided that $f'(\overline{u})F(\overline{u})\le A$. 
Then, with the help of the cut-off technique, we can 
construct a supersolution by using \eqref{eq:1.24} and the solution of \eqref{eq:1.25}. 
Nonexistence of solutions for problem~\eqref{eq:1.1} is also proved by using \eqref{eq:1.24} 
with the help of the nonexistence results for \eqref{eq:1.25}. 
However, the transformation~\eqref{eq:1.24} is not useful for the case $A=1$ 
since the case $A=1$ includes exponential nonlinearity $f(u)=e^{u}$, which is essentially different from power type nonlinearity. 
For the case $A=1$, we use the transformation 
\begin{equation} 
\label{eq:1.26} 
v(x,t):=\log F(u(x,t))^{-1}, 
\end{equation} 
instead of \eqref{eq:1.24}. 
Under this transformation, the existence problem for \eqref{eq:1.1} with general $f$ can be reduced to that of a heat equation with exponential nonlinearity.  
\vspace{5pt} 

The rest of this paper is organized as follows: 
In Section~\ref{section:2}, we give some preliminary results. 
In particular, we recall the existence and nonexistence results for a heat equation with power type nonlinearity and exponential nonlinearity. 
In Section~\ref{section:3}, we consider local 
in time existence of solutions for problem~\eqref{eq:1.1} with the aid of \eqref{eq:1.24} and \eqref{eq:1.26}, 
and prove Theorems~\ref{Theorem:1.1}. 
In Section~\ref{section:4}, we discuss nonexistence of solutions for problem~\eqref{eq:1.1}, 
and prove Theorems~\ref{Theorem:1.3} and \ref{Theorem:1.4}. 
In Section~\ref{section:5}, we apply our main theorems to several examples of nonlinear heat equations. 

\section{Preliminaries} 
\label{section:2} 
In this section we recall some properties of  
uniformly local $L^p$ spaces and the existence result of solutions for problem~\eqref{eq:1.1}. 
Furthermore, in Propositions~\ref{Proposition:2.2}--\ref{Proposition:2.5}, 
we discuss the existence result of solutions for problem~\eqref{eq:1.1} with typical examples of $f$. 
In particular, we discuss the cases $f(u)=u^p$ and $f(u)=e^u$. 
\vspace{5pt} 

We first recall two lemmas on properties of uniformly local $L^p$ spaces. 
For partial differential equations in the uniformly local Lebesgue spaces, 
see for example \cite{ABCD}, \cite{IS} and \cite{MT}. 
Lemma~\ref{Lemma:2.1} gives the smoothing effect of the heat semigroup in $L^p_{ul,\rho}(\mathbb{R}^N)$ 
(see \cite[Proposition~2.1]{ABCD} and \cite[Corollary~3.1]{MT}). 
In the following, for any set $X$ and maps $a=a(x)$ and $b=b(x)$ from $X$ to $[0,\infty)$, we say 
$$
a(x)\lesssim b(x) 
\quad\mbox{for all}\,\,\, x \in X, 
$$
if there exists a positive constant $C$ such that $a(x)\le Cb(x)$ 
for all $x\in X$. 
\begin{Lemma} 
\label{Lemma:2.1}
Let $1\le p \le q \le \infty$.
Then there holds 
$$
\|e^{t\Delta}u\|_{L^q_{ul,\rho}(\mathbb{R}^N)}\lesssim 
\left(
\rho^{-N\left(\frac{1}{p}-\frac{1}{q}\right)} 
+t^{-\frac{N}{2}\left(\frac{1}{p}-\frac{1}{q}\right)}\right)\|u\|_{L^p_{ul,\rho}(\mathbb{R}^N)}
$$
for all $t>0$, $\rho>0$ and $u\in L^{p}_{ul,\rho}(\mathbb{R}^N)$.
\end{Lemma}
The following lemma gives 
basic properties of $\mathcal{L}^p_{ul,\rho}(\mathbb{R}^N)$. 
See \cite[Proposition~2.2]{MT}. 
\begin{Lemma} 
\label{Lemma:2.2}
Let $1\le p<\infty$. 
The following assertions are  equivalent: 
$$
\begin{array}{cl} 
{\rm (i)} &u\in {\mathcal{L}^p_{ul,\rho}(\mathbb{R}^N)}.
\vspace{3pt} \\
{\rm (ii)} &\displaystyle\lim_{|y|\to 0}\|u(\cdot +y)-u(\cdot)\|_{L^p_{ul,\rho}(\mathbb{R}^N)}=0.
\vspace{3pt} \\
{\rm (iii)} &\displaystyle\lim_{t\to 0}\|e^{t\Delta}u-u\|_{L^p_{ul,\rho}(\mathbb{R}^N)}=0.
\end{array} 
$$ 
\end{Lemma} 
\begin{Remark} 
\label{Remark:2.1} 
$u\in {\cal L}^p_{ul,\rho}(\mathbb{R}^N)$ is equivalent to $|u|^p\in {\cal L}_{ul,\rho}^1(\mathbb{R}^N)$. 
\end{Remark} 
It is pointed out in \cite{MT} that 
the characterization (iii) plays an important role to treat the initial value problem in 
$\mathcal{L}^{p}_{ul,\rho}(\mathbb{R}^N)$ 
for a nonlinear heat equation. 
\vspace{5pt} 

We next recall one proposition on existence of solutions for problem~\eqref{eq:1.1}. 
Proposition~\ref{Proposition:2.1} implies that, if there exists a supersolution for \eqref{eq:1.1}, 
we can find a solution of \eqref{eq:1.1} below the supersolution via monotone methods. 
See for example \cite{IKS}, \cite{RS} and \cite{S}. 
\begin{Proposition} 
\label{Proposition:2.1} 
Let $f\in C^1([0,\infty))$ satisfy \eqref{eq:1.2}. 
Let $u_0\in L_{ul,\rho}^1({\mathbb{R}}^N)$ be a nonnegative function. 
Assume that there exists 
$\overline{u}\in C^{2,1}(\mathbb{R}^N\times (0,T))$ satisfying 
\begin{equation} 
\label{eq:2.1} 
\overline{u}(x,t)\ge (e^{t\Delta}u_0)(x)+\int_0^t [e^{(t-s)\Delta}f(\overline{u}(\cdot, s))](x) \, ds 
\end{equation} 
for almost every $(x,t)\in \mathbb{R}^N\times (0,T)$. 
Then there exists a solution $u\in C^{2,1}({\mathbb R}^N\times (0,T))$ of the integral equation 
\begin{equation} 
\notag 
u(x,t)=(e^{t\Delta}u_0)(x)+\int_0^t [e^{(t-s)\Delta}f(u(\cdot, s))](x)\, ds \quad\mbox{in}\,\,\, \mathbb{R}^N\times (0,T),  
\end{equation} 
which satisfies $\partial_t u = \Delta u + f(u)$ in $\mathbb{R}^N\times (0,T)$. 
Furthermore, there holds $0\le u(x,t)\le \overline{u}(x,t)$ in ${\mathbb R}^N\times (0,T)$. 
\end{Proposition} 
\noindent 
{\bf Proof.} 
For $n\ge 2$, define the function $u_n$ by  
$$
u_n(x,t):=(e^{t\Delta}u_{0})(x)+\int_0^t [e^{(t-s)\Delta}f(u_{n-1}(s))](x)\, ds,
$$
where $u_1:=0$. 
Since $f=f(u)$ is an increasing function with respect to $u$ by \eqref{eq:1.2} and $\overline{u}$ is a supersolution 
in the sense of \eqref{eq:2.1}, 
if $ \overline{u}\ge u_{n-1}(x,t)$, then we have 
\begin{align*} 
\overline{u}(x,t) 
&\ge (e^{t\Delta}u_0)(x)+\int_0^t [e^{(t-s)\Delta}f(\overline{u}(s))](x)\, ds \\ 
&\ge (e^{t\Delta}u_0)(x)+\int_0^t [e^{(t-s)\Delta}f(u_{n-1}(s))](x)\, ds=u_n(x,t). 
\end{align*} 
Since $\overline{u}(x,t)\ge 0=u_1(x,t)$, by induction we have $\overline{u}(x,t)\ge u_n(x,t)$ for all $n\in \mathbb{N}$. 
By the definition of $u_n$ with $n=2$ we first obtain 
$$
u_2(x,t)=(e^{t\Delta}u_{0})(x)+\int_0^t [e^{(t-s)\Delta}f(u_{1}(s))](x)\, ds 
\ge 0=u_1(x,t). 
$$
Then, since $f=f(u)$ is an increasing function with respect to $u$, we have 
$$ 
u_3=e^{t\Delta}u_{0}+\int_0^t e^{(t-s)\Delta}f(u_{2}(s))\, ds 
\ge 
e^{t\Delta}u_{0}+\int_0^t e^{(t-s)\Delta}f(u_{1}(s))\, ds 
=u_2. 
$$
Repeating the above argument, we have 
\begin{equation} 
\label{eq:2.2} 
\overline{u}(x,t)\ge u_{n+1}(x,t)\ge u_n(x,t) \ge 0 \quad\mbox{in}\,\,\, \mathbb{R}^N\times (0,T) 
\end{equation} 
for all $n\in \mathbb{N}$. 
Then we can define the limit function $u$ by 
$$
u(x,t):=\lim_{n\to\infty}u_n(x,t), 
$$
and by the monotone convergence theorem and \eqref{eq:1.2} we see that this function $u$ gives a solution to the desired integral equation. 
By \eqref{eq:2.2} and the monotonicity of $u_n$ we have $u(x,t)\le \overline{u}(x,t)$. 
Then we can apply the standard regularity theory for parabolic equations and obtain $u\in C^{2,1}(\mathbb{R}^N\times (0,T))$, 
so $u$ satisfies the equation in the classical sense. 
Thus we complete the proof of Proposition~\ref{Proposition:2.1}. 
\hfill \qed 

\vspace{5pt} 

We give one lemma on sufficient conditions that equation~\eqref{eq:1.1} and related inequalities can be rewritten by the integral form. 
%
\begin{Lemma} 
\label{Lemma:2.3} 
Let $f\in C^1([0,\infty))$ satisfy \eqref{eq:1.2}.
Let $T>0$ and $u\in C^{2,1}(\mathbb{R}^N\times (0,T))$ satisfy 
\begin{equation} 
\label{eq:2.3} 
\lim_{t\to 0} 
\left\| 
\int_0^t e^{(t-s)\Delta}f(u(\cdot, s))\, ds 
\right\|_{L_{ul,\rho}^1(\mathbb{R}^N)} 
=0. 
\end{equation} 
\begin{itemize} 
\item[\rm (i)] 
Assume that $u$ satisfies $\partial_t u=\Delta u+f(u)$ in $\mathbb{R}^N\times (0,T)$. 
If $u_0\in L_{ul,\rho}^1(\mathbb{R}^N)$ and there holds 
\begin{equation} 
\label{eq:2.4} 
\lim_{t\to 0}\|u(t)-e^{t\Delta}u_0\|_{L_{ul,\rho}^1(\mathbb{R}^N)}=0,
\end{equation} 
then $u$ satisfies 
$$
u(x,t)=(e^{t\Delta}u_0)(x)+\int_0^t [e^{(t-s)\Delta}f(u(\cdot, s))](x)\, ds 
\quad\mbox{in}\,\,\, \mathbb{R}^N\times (0,T). 
$$ 
\item[\rm (ii)] 
Assume that $u$ satisfies $\partial_t u\ge \Delta u+f(u)$ in $\mathbb{R}^N\times (0,T)$. 
If $u_0\in L_{ul,\rho}^1(\mathbb{R}^N)$ and there holds either \eqref{eq:2.4}
or 
\begin{equation} 
\label{eq:2.5} 
u(x,t)\ge (e^{t\Delta}u_0)(x) \quad\mbox{in}\,\,\, \mathbb{R}^N\times (0,T), 
\end{equation} 
then $u$ satisfies 
$$
u(x,t)\ge (e^{t\Delta}u_0)(x)+\int_0^t [e^{(t-s)\Delta}f(u(\cdot, s))](x)\, ds 
\quad\mbox{in}\,\,\, \mathbb{R}^N\times (0,T). 
$$ 
\item[\rm (iii)] 
Assume that $u$ satisfies $\partial_t u\le \Delta u+f(u)$ in $\mathbb{R}^N\times (0,T)$. 
If $u_0\in L_{ul,\rho}^1(\mathbb{R}^N)$ and there holds either \eqref{eq:2.4}
or that 
there exists a function $r\in C([0,T); L_{ul,\rho}^1(\mathbb{R}^N))$ such that 
$$ 
  \|r(\cdot,t)\|_{L_{ul,\rho}^1(\mathbb{R}^N)}\to 0 \quad \mbox{as} \,\,\, t\to 0
$$ 
and 
\begin{equation} 
\label{eq:2.6} 
u(x,t)\le (e^{t\Delta}u_0)(x)+r(x,t) \quad\mbox{in}\,\,\, \mathbb{R}^N\times (0,T), 
\end{equation} 
then $u$ satisfies 
$$
u(x,t)\le (e^{t\Delta}u_0)(x)+\int_0^t [e^{(t-s)\Delta}f(u(\cdot, s))](x)\, ds 
\quad\mbox{in}\,\,\, \mathbb{R}^N\times (0,T). 
$$ 
\end{itemize} 
\end{Lemma} 
\noindent 
{\bf Proof.} 
We first prove assertion~(i). 
Since $u$ satisfies $\partial_t u=\Delta u+f(u)$ in $\mathbb{R}^N\times (0,T)$, for any sufficiently small $\tau>0$, we have 
\begin{equation} 
\label{eq:2.7} 
u(x,t)=(e^{(t-\tau)\Delta}u(\tau))(x)+\int_{\tau}^t [e^{(t-s)\Delta}f(u(s))](x)\, ds \quad\mbox{in}\,\,\, \mathbb{R}^N\times (\tau, T). 
\end{equation} 
Since $e^{(t-\tau)\Delta}u(\tau)=e^{(t-\tau)\Delta}(u(\tau)-e^{\tau\Delta}u_0)+e^{t\Delta}u_0$, 
by Lemma~\ref{Lemma:2.1} and \eqref{eq:2.4} we have 
\begin{equation} 
\label{eq:2.8} 
\begin{split} 
& \lim_{\tau\to 0}\|e^{(t-\tau)\Delta}u(\tau)-e^{t\Delta}u_0\|_{L^\infty (\mathbb{R}^N)} \\ 
&\qquad 
\lesssim \lim_{\tau\to 0}(\rho^{-N}+(t-\tau)^{-\frac{N}{2}})\|u(\tau)-e^{\tau\Delta}u_0\|_{L_{ul,\rho}^1(\mathbb{R}^N)} =0. 
\end{split} 
\end{equation} 
On the other hand, since 
$$
\int_0^t e^{(t-s)\Delta}f(u(s))\, ds-\int_{\tau}^t e^{(t-s)\Delta}f(u(s))\, ds
=e^{(t-\tau)\Delta} \int_0^\tau e^{(\tau-s)\Delta}f(u(s))\, ds, 
$$
by Lemma~\ref{Lemma:2.1} and \eqref{eq:2.3} we obtain 
\begin{align*} 
&\lim_{\tau \to 0} \left\| 
\int_0^t e^{(t-s)\Delta}f(u(s))\, ds-\int_{\tau}^t e^{(t-s)\Delta}f(u(s))\, ds
\right\|_{L^\infty(\mathbb{R}^N)} \\ 
&\qquad \lesssim \lim_{\tau \to 0} (\rho^{-N}+(t-\tau)^{-\frac{N}{2}})\left\| 
\int_0^\tau e^{(\tau-s)\Delta}f(u(s))\, ds 
\right\|_{L_{ul,\rho}^1(\mathbb{R}^N)} =0. 
\end{align*} 
This together with \eqref{eq:2.7} and \eqref{eq:2.8} proves assertion~(i). 

Next we prove assertion~(ii). 
If $u_0$ satisfies \eqref{eq:2.4}, then we can prove assertion~(ii) as in the above argument. 
If \eqref{eq:2.5} is satisfied, we have 
$$
u(x,t)\ge e^{(t-\tau)\Delta}u(\tau)+\int_{\tau}^t e^{(t-s)\Delta}f(u(s))\, ds 
\ge e^{t\Delta}u_0+\int_{\tau}^t e^{(t-s)\Delta}f(u(s))\, ds  
$$ 
in $\mathbb{R}^N\times (\tau, T)$, where $\tau>0$ is sufficiently small. 
Then, as in the above argument, we have the convergence of the Duhamel term 
and assertion~(ii) is proved. 

Finally, we prove assertion~(iii). 
We only consider the case \eqref{eq:2.6}. 
Assuming \eqref{eq:2.6}, we have 
$$
u(t)\le e^{(t-\tau)\Delta}u(\tau)+\int_{\tau}^t e^{(t-s)\Delta}f(u(s))\, ds 
\le e^{t\Delta}u_0+e^{(t-\tau)\Delta}r(\tau)+\int_{\tau}^t e^{(t-s)\Delta}f(u(s))\, ds  
$$ 
in $\mathbb{R}^N\times (\tau, T)$, where $\tau>0$ is sufficiently small. 
Since $\|r(t)\|_{L_{ul,\rho}^1(\mathbb{R}^N)}\to 0$ as $t\to 0$, by Lemma~\ref{Lemma:2.1} we have 
$$
\|e^{(t-\tau)\Delta}r(\tau)\|_{L^\infty(\mathbb{R}^N)}\lesssim (\rho^{-N}+(t-\tau)^{-\frac{N}{2}})\|r(\tau)\|_{L_{ul,\rho}^1(\mathbb{R}^N)}\to 0 
$$
as $\tau\to 0$. 
Since the convergence of the Duhamel term can be proved as in the above argument, 
we can prove assertion~(iii). 
Thus we complete the proof of Lemma~\ref{Lemma:2.3}. 
\hfill \qed 

\vspace{5pt} 

Now we recall the existence result for the heat equation with power type nonlinearity 
\begin{equation} 
\label{eq:2.9} 
\left\{ 
\begin{array}{ll} 
\partial_t u=\Delta u+|u|^{p-1}u, &x\in {\mathbb{R}}^N, \,\,\, t>0, \vspace{3pt} \\ 
u(x,0)=u_0(x), &x\in {\mathbb{R}}^N, 
\end{array} 
\right. 
\end{equation} 
where $p>1$. 
In particular, we consider the case $u_0$ belongs to a uniformly local $L^r$ space and 
study local in time existence of solutions for problem~\eqref{eq:2.9} in suitable functional spaces. 
We state the existence results for the subcritical case and 
the critical case, respectively. 
For existence of classical solutions of \eqref{eq:2.9}, we study the integral equation 
\begin{equation} 
\label{eq:2.10} 
u(t)=e^{t\Delta}u_0+\int_0^t e^{(t-s)\Delta}|u(s)|^{p-1}u(s)\, ds 
\end{equation} 
in the uniformly local $L^r$ spaces. 
\begin{Proposition}[Subcritical case] 
\label{Proposition:2.2} 
Let $N\ge 1$ and $p>1$. 
Assume $r\ge 1$ and $r>\frac{N}{2}(p-1)$. 
Given any $u_0\in {L}_{ul,\rho}^r({\mathbb{R}}^N)$, 
there exist $T>0$ and at least one classical solution 
$u\in C((0,T); {L}_{ul,\rho}^r({\mathbb{R}}^N))\cap L^{\infty}_{loc}((0,T); L^{\infty}(\mathbb{R}^N))\cap C^{2,1}(\mathbb{R}^N\times (0,T))$ 
satisfying \eqref{eq:2.10} and 
$$
\lim_{t\to 0}\|u(t)-e^{t\Delta}u_0\|_{L^r_{ul,\rho}(\mathbb{R}^N)}=0. 
$$
The solution $u$ satisfies $\|u(t)-u_0\|_{L^r_{ul,\rho}(\mathbb{R}^N)}\to 0$ as $t\to 0$, 
provided that $u_0\in {\mathcal L}^r_{ul,\rho}(\mathbb{R}^N)$. 
Furthermore, the maximal existence time $T$ can be estimated to satisfy 
\begin{equation} 
\label{eq:2.11} 
 T^{\frac{N}{2p}} \rho^{-\frac{N}{p}} 
 + \|u_0\|_{L^{r}_{ul,\rho}(\mathbb{R}^N)}^r
 \left( 
 T^{\frac{r}{p-1}-\frac{N}{2}\cdot\frac{p-1}{p}}\rho^{-\frac{N}{p}}+T^{\frac{r}{p-1}-\frac{N}{2}}
 \right) 
 \ge \gamma, 
\end{equation}
where $\gamma$ is a positive constant depending only on $N$, $p$ and $r$. 
\end{Proposition} 
\begin{Proposition}[Critical case] 
\label{Proposition:2.3}
Let $N\ge 1$ and assume $r=\frac{N}{2}(p-1)>1$.
Given any $u_0\in \mathcal{L}_{ul,\rho}^r({\mathbb{R}}^N)$, 
there exist $T>0$ and at least one classical solution 
$u\in C([0,T);\mathcal{L}_{ul,\rho}^r({\mathbb{R}}^N))\cap L^{\infty}_{loc}((0,T);L^{\infty}(\mathbb{R}^N)) 
\cap C^{2,1}(\mathbb{R}^N\times (0, T))$ satisfying \eqref{eq:2.10}.  
\end{Proposition}
One can prove Propositions~\ref{Proposition:2.2} and \ref{Proposition:2.3} 
by applying the arguments in \cite{GV} and \cite{W1} 
with a slight modification. 
See Appendix~\ref{section:A}. 
\begin{Remark} 
\label{Remark:2.2} 
{\rm (i)} 
Under the assumptions of Proposition~{\rm\ref{Proposition:2.2}}, there exist constants $M>0$ and $T_0\in (0,T)$ such that 
$$
t^\sigma \|u(t)\|_{L_{ul,\rho}^{pr}(\mathbb{R}^N)}\le M 
$$
for all $t\in (0,T_0)$, 
where $\sigma=\frac{N}{2}(\frac{1}{r}-\frac{1}{pr})$. 
See the proof of Proposition~{\rm\ref{Proposition:2.2}} in Appendix~\ref{section:A}. 
\vspace{3pt} 
\newline 
{\rm (ii)} 
Under the assumptions of Proposition~{\rm\ref{Proposition:2.3}}, 
it holds 
$$
\lim_{t\to 0}\, t^\sigma \|u(t)\|_{L_{ul,\rho}^q (\mathbb{R}^N)}=0, 
$$
where $\max\{p,r\}<q<pr$ and $\sigma=\frac{N}{2}(\frac{1}{r}-\frac{1}{q})$. 
See the proof of Proposition~{\rm\ref{Proposition:2.3}} in Appendix~\ref{section:A}. 
\end{Remark} 
%
%
\vspace{5pt} 

We next discuss existence of solutions for a heat equation with exponential nonlinearity 
\begin{equation} 
\label{eq:2.12} 
\left\{ 
\begin{array}{ll} 
\partial_t u=\Delta u+e^u, &x\in {\mathbb{R}}^N, \,\,\, t>0, \vspace{3pt} \\ 
u(x,0)=u_0(x), &x\in {\mathbb{R}}^N, 
\end{array} 
\right. 
\end{equation} 
where $u_0$ satisfies 
\begin{equation} 
\label{eq:2.13} 
u_0(x)\ge -Ce^{|x|^{2-\epsilon}} \quad\mbox{in}\,\,\, \mathbb{R}^N 
\end{equation} 
for some $C>0$ and $\epsilon\in (0,1)$. 
The problem on existence and nonexistence of solutions for \eqref{eq:1.1} with $A=1$ can be reduced to that of \eqref{eq:2.12} 
via the transformation~\eqref{eq:1.26} 
under the condition \eqref{eq:1.9} or \eqref{eq:1.20}. 

We first prepare one basic lemma on the relationship between 
the heat semigroup and convex and concave functions. 
Lemma~\ref{Lemma:2.4} directly follows from the Jensen inequality. 
For the proof of Lemma~\ref{Lemma:2.4}, we refer to \cite[Lemma~5.1]{W1}. 
\begin{Lemma}\label{Lemma:2.4}
Let $\phi\in L^1_{ul,\rho}(\mathbb{R}^N)$
and $J$ be a function from $[0, \infty)$ to itself. 
If $J$ is convex,
then
 there holds
$$
 J\left(e^{t\Delta}\phi(x)\right)\le [e^{t\Delta}J(\phi)](x)\ \  \mbox{in}\ \mathbb{R}^N\times (0,\infty).
$$
On the other hand, if $J$ is concave,
then
 there holds
$$
 J\left(e^{t\Delta}\phi(x)\right)\ge [e^{t\Delta}J(\phi)](x) \ \  \mbox{in}\ \mathbb{R}^N\times (0,\infty).
$$
\end{Lemma} 
%

%
%
\begin{Proposition} 
\label{Proposition:2.4}
Let $N\ge 1$ and $r\ge N/2$. 
For any (possibly sign changing) initial data $u_0$ satisfying \eqref{eq:2.13} and 
$$
e^{ru_0}\in L^1_{ul,\rho}({\mathbb{R}}^N) \quad\mbox{if} \,\,\, r>\frac{N}{2}, \qquad 
e^{ru_0}\in {\cal L}^1_{ul,\rho}(\mathbb{R}^N) \quad\mbox{if} \,\,\, r=\frac{N}{2}, 
$$
there exists at least one classical solution $u$ of \eqref{eq:2.12} satisfying 
$$
u(t)=e^{t\Delta}u_0+\int_0^t e^{(t-s)\Delta}e^{u(s)}\, ds  
\quad\mbox{and}\quad 
\lim_{t\to 0} \|u(t)-e^{t\Delta}u_0\|_{L^{\infty}(\mathbb{R}^N)}=0. 
$$
If $r>\frac{N}{2}$, then the existence time $T$ can be taken to satisfy
\begin{equation}\label{eq:2.14} 
 T^{\frac{\epsilon N}{2}}\rho^{-\epsilon N} 
 + \gamma_\epsilon \|e^{ru_0}\|_{L_{ul,\rho}^1(\mathbb{R}^N)} 
 \left(T^{r-\frac{N}{2}(1-\epsilon)}\rho^{-\epsilon N}+T^{r-\frac{N}{2}}\right) 
 \ge 
 \gamma  
 \end{equation}
for any sufficiently small $\epsilon>0$, 
where $\gamma$ depends only on $N$ and $r$, 
and $\gamma_\epsilon$ is a positive constant depending only on $N$, $r$ and $\epsilon$ 
satisfying $\gamma_\epsilon\to \infty $ as $\epsilon\to 0$. 
\end{Proposition} 

\noindent 
{\bf Proof.} 
Let $\epsilon\in (0,1/3)$ and $v_0:= e^{\frac{\epsilon}{1-\epsilon}u_0}$. 
Put $X_r:=L^{\frac{(1-\epsilon) r}{\epsilon}}_{ul,\rho}(\mathbb{R}^N)$ if $r>N/2$ 
and $X_r:=\mathcal{L}^{\frac{(1-\epsilon) r}{\epsilon}}_{ul,\rho}(\mathbb{R}^N)$ if $r=N/2$. 
By the assumption on $u_0$ we have $v_0\in X_r$ for $r\ge N/2$. 
Consider 
\begin{equation} 
\label{eq:2.15} 
\partial_t v=\Delta v+\frac{\epsilon }{1-\epsilon}v^{1/\epsilon}, \quad x\in \mathbb{R}^N,\,\,\, t>0, \qquad 
v(x,0)=v_0(x), \quad x\in \mathbb{R}^N. 
\end{equation} 
Since $\frac{(1-\epsilon) r}{\epsilon}\ge \frac{1-\epsilon}{\epsilon}\cdot \frac{N}{2}=\frac{N}{2}(\frac{1}{\epsilon}-1)>1$, 
Propositions~\ref{Proposition:2.2} and \ref{Proposition:2.3} imply existence of a classical nonnegative solution 
$v$ of \eqref{eq:2.15} in $C((0,T]; X_r)$
if $r>N/2$ and in $C([0,T]; X_r)$ 
if $r=N/2$,
for some $T>0$.
Now we apply the Cole-Hopf transformation
$$
 \overline{u}(x,t):=\frac{1-\epsilon}{\epsilon} \log v(x,t).
$$
Then a simple calculation shows that $v^{\frac{1-\epsilon}{\epsilon}}=e^{\overline{u}}$ and 
\begin{equation}
\label{eq:2.16}
\partial_t \overline{u}-\Delta \overline{u} =\frac{1-\epsilon}{\epsilon} 
\left( 
\frac{\partial_t v-\Delta v}{v}+\frac{|\nabla v|^2}{v^2} 
\right) 
\ge e^{\overline u} 
\quad\mbox{in}\,\,\, \mathbb{R}^N\times (0,T). 
\end{equation}
Therefore $\overline u$ is a supersolution of \eqref{eq:2.12}. 
In order to rewrite \eqref{eq:2.16} by the integral form, 
we check the assumptions of Lemma~\ref{Lemma:2.3}~(ii). 
By \eqref{eq:2.10} we have $v(t)\ge e^{t\Delta}v_0$. 
Since $\log s$ is monotonically increasing and concave with respect to $s$, by Lemma~\ref{Lemma:2.4} 
we have 
\begin{equation} 
\label{eq:2.17} 
\overline{u}(t)-e^{t\Delta}\overline{u}(0) 
=\frac{1-\epsilon}{\epsilon} (\log v(t)-e^{t\Delta}(\log v_0))  
\ge \frac{1-\epsilon}{\epsilon}( \log v(t)- \log (e^{t\Delta}v_0)) \ge 0,
\end{equation} 
thus $\overline{u}$ satisfies condition \eqref{eq:2.5}. 
We now check that $\overline{u}$ satisfies condition \eqref{eq:2.3} with $f(u)=e^u$. 
By the definition of $\overline{u}$, 
it suffices to prove
\begin{equation}\label{eq:2.18}
 \int_0^t \| e^{(t-s)\Delta}v(s)^{\frac{1-\epsilon}{\epsilon}}\|_{L^{\infty}(\mathbb{R}^N)}\, ds \to 0\ \ \mbox{as}\ \ t\to 0. 
\end{equation}

\noindent 
\underline{\bf Case $r>N/2$} 
\vspace{3pt} 
\newline 
We remark that $v_0\in L_{ul,\rho}^{\frac{1-\epsilon}{\epsilon} r}(\mathbb{R}^N)$. 
If $r>N/2$, then there exist constants $M>0$ and $T_0\in (0,T)$ such that
$t^\sigma \|v(t)\|_{L^{\frac{1-\epsilon}{\epsilon^2}r}_{ul,\rho}(\mathbb{R}^N)}\le M$ for $0<t<T_0$, 
where $\sigma=\frac{N}{2}(\frac{1}{\frac{1-\epsilon}{\epsilon} r}-\frac{1}{\frac{1-\epsilon}{\epsilon^2}r})=\frac{\epsilon N}{2r}$. 
See Remark~\ref{Remark:2.2}~(i). 
Thus, by Lemma~\ref{Lemma:2.1} we have 
\begin{equation} 
\label{eq:2.19} 
\begin{split} 
 \int_0^t \| e^{(t-s)\Delta}v(s)^{\frac{1-\epsilon}{\epsilon}}\|_{L^{\infty}(\mathbb{R}^N)}\, ds 
 & \lesssim \int_0^t \left(\rho^{-2\sigma}+(t-s)^{-\sigma}\right)\|v(s)\|^{\frac{1-\epsilon}{\epsilon}}_{L^{\frac{1-\epsilon}{\epsilon^2}r}_{ul,\rho}(\mathbb{R}^N)}\, ds \\ 
 & \le \int_0^t  \left(\rho^{-2\sigma}+(t-s)^{-\sigma}\right) s^{-\frac{1-\epsilon}{\epsilon} \sigma}\, ds \cdot M^\frac{1-\epsilon}{\epsilon} \\ 
 & \lesssim ({t^{1-\frac{1-\epsilon}{\epsilon}\sigma}}\rho^{-2\sigma}+t^{1-\frac{1-\epsilon}{\epsilon} \sigma-\sigma})M^{\frac{1-\epsilon}{\epsilon}}\to 0
\end{split} 
\end{equation}  
as $t\to 0$, since $\sigma<\frac{1-\epsilon}{\epsilon} \sigma=\frac{(1-\epsilon)N}{2r}=\frac{N}{2r}-\sigma<1$ 
and $1-\-\frac{1-\epsilon}{\epsilon}\sigma-\sigma=1-\frac{\sigma}{\epsilon}=1-\frac{N}{2r}>0$. 
\vspace{5pt} 
\newline 
\underline{\bf Case $r=N/2$} 
\vspace{3pt} 
\newline 
In the case $r=N/2$, by Remark~\ref{Remark:2.2}~(ii) we have $t^{\sigma}\|v(t)\|_{L^{\alpha}_{ul,\rho}(\mathbb{R}^N)}\to 0$ as $t\to 0$, 
where $\max\left\{\frac{1-\epsilon}{\epsilon} r, \frac{1}{\epsilon}\right\}<\alpha<\frac{1-\epsilon}{\epsilon^2}r$ 
and $\sigma=\frac{N}{2}(\frac{\epsilon}{1-\epsilon}\cdot \frac{1}{r}-\frac{1}{\alpha})$.
Therefore, by Lemma~\ref{Lemma:2.1} we have 
\begin{equation} 
\label{eq:2.20} 
\begin{split} 
 & 
 \int_0^t \| e^{(t-s)\Delta}v(s)^{\frac{1-\epsilon}{\epsilon}}\|_{L^{\infty}(\mathbb{R}^N)}\, ds \\ 
 &
 \quad \lesssim \int_0^t \left(\rho^{-\frac{1-\epsilon}{\epsilon}\cdot\frac{N}{\alpha}}+(t-s)^{-\frac{1-\epsilon}{\epsilon}\cdot \frac{N}{2\alpha}}\right) 
 s^{-\frac{1-\epsilon}{\epsilon}\cdot\sigma}\, ds 
 \cdot \sup_{0<s<t}s^{\frac{1-\epsilon}{\epsilon}\sigma}\|v(s)\|_{L_{ul,\rho}^{\alpha}(\mathbb{R}^N)}^{\frac{1-\epsilon}{\epsilon}}  
 \\
 &
 \quad \lesssim 
 (t^{1-\frac{1-\epsilon}{\epsilon}\sigma}\rho^{-\frac{1-\epsilon}{\epsilon}\cdot \frac{N}{\alpha}}+1)\cdot 
 \sup_{0<s<t}s^{\frac{1-\epsilon}{\epsilon}\sigma}\|v(s)\|_{L_{ul,\rho}^{\alpha}(\mathbb{R}^N)}^{\frac{1-\epsilon}{\epsilon}} \to 0
\end{split} 
\end{equation} 
as $t\to 0$, since $-\frac{1-\epsilon}{\epsilon}\cdot \frac{N}{2\alpha}>-1$, $-\frac{1-\epsilon}{\epsilon} \sigma>-1$ 
and $-\frac{1-\epsilon}{\epsilon}\cdot \frac{N}{2\alpha}-\frac{1-\epsilon}{\epsilon}\sigma+1=0$. 
\vspace{5pt} 

Hence \eqref{eq:2.16} can be written by the integral form
by Lemma~\ref{Lemma:2.3}~(ii) with the aid of \eqref{eq:2.17} and \eqref{eq:2.18}. 
Then, by Proposition~\ref{Proposition:2.1} we obtain a classical solution $u$ of
\begin{equation} 
\label{eq:2.21} 
u(x,t)=(e^{t\Delta}u_0)(x)+\int_0^t \left[e^{(t-s)\Delta}e^{u(s)}\right](x)\, ds \quad\mbox{in}\,\,\, \mathbb{R}^N\times (0,T), 
\end{equation} 
satisfying $u(x,t)\le \overline{u}(x,t)$. 
We now prove the convergence of a solution to the initial data. 
Since $u\le \overline{u}$, by \eqref{eq:2.21} we have 
$$
 \|u(t)-e^{t\Delta}u_0\|_{L^{\infty}(\mathbb{R}^N)} 
 \le \left\|\int_0^t e^{(t-s)\Delta}e^{\overline u (s)}\, ds\right\|_{L^{\infty}(\mathbb{R}^N)} 
 \le  \int_0^t \| e^{(t-s)\Delta}v(s)^{\frac{1-\epsilon}{\epsilon}}\|_{L^{\infty}(\mathbb{R}^N)}\, ds, 
$$
and obtain the convergence of $u$ to the initial data by \eqref{eq:2.19} and \eqref{eq:2.20}. 

We finally study the estimate of the existence time $T$ for the case $r>N/2$, 
and prove \eqref{eq:2.14}. 
Let $r>N/2$ and $\epsilon\in (0,1/3)$. 
One can apply the same argument as in Appendix~\ref{section:A} 
with $p$, $r$, $\alpha$ ,$u^p$ 
replaced by $1/\epsilon$, $\frac{1-\epsilon}{\epsilon}r$, $\frac{\epsilon N}{2r}$, $\frac{\epsilon}{1-\epsilon}v^{1/\epsilon}$, respectively. 
Then, being careful with the changes of the constants, we see that 
$$
T^{r-\frac{N}{2}+\epsilon\frac{N}{2}}\rho^{-\epsilon N}+T^{r-\frac{N}{2}}\ge 
\frac{\gamma\epsilon^{2r}}{\|v_0\|_{L_{ul,\rho}^{\frac{1-\epsilon}{\epsilon}r}(\mathbb{R}^N)}^{\frac{1-\epsilon}{\epsilon}r}}, 
$$ 
which proves \eqref{eq:2.14} with the help of the definition of $v_0$. 
Thus we complete the proof of Proposition~\ref{Proposition:2.4}. 
\hfill 
\qed 
	
\vspace{5pt} 

The following result states nonexistence of solutions for 
exponential nonlinear heat equation including problem~\eqref{eq:2.12}, 
which shows the optimality of the condition for the integrability of the initial data in Proposition~\ref{Proposition:2.4}. 
Proposition~\ref{Proposition:2.5} is also available for rapidly increasing nonlinearity such as $f(u)=e^{u^2}$, 
and is the key assertion for nonexistence of solutions for problem~\eqref{eq:1.1} 
even if \eqref{eq:1.20} is violated. 
\begin{Proposition} 
\label{Proposition:2.5}
Let $r\in (0,N/2)$. 
Let $g$ be a convex function in $(s_0, \infty)$ for some $s_0>0$. 
Assume that $g$ satisfies 
$g(s)\to \infty$ as $s\to\infty$, 
$g'(s)>0$ for all $s>0$ and 
\begin{equation} 
\label{eq:2.22} 
\lim_{s\to\infty}\frac{g''(s)}{(g'(s))^2} = 0. 
\end{equation} 
Then there exists $u_0\ge 0$ satisfying 
$$
G(u_0)^{-r} \in {L}^1_{ul,\rho}(\mathbb{R}^N) 
\quad\mbox{with} \quad 
G(s):=\int_s^\infty \frac{du}{e^{g(u)}} 
$$
such that, for every $T>0$, 
there is no nonnegative solution $u\in C^{2,1}(\mathbb{R}^N \times (0,T))$
of
$$
 u(x,t)=e^{t\Delta}u_0 +\int_0^t e^{(t-s)\Delta}e^{g(u(s))}ds \quad\mbox{in}\,\,\, \mathbb{R}^N\times (0,T).
$$
In particular, there is no nonnegative classical solution $u\in C^{2,1}(\mathbb{R}^N \times (0,T))$ of 
\begin{equation} 
\label{eq:2.23} 
	\partial_t u = \Delta u + e^{g(u)} 
	\quad\mbox{in} \,\,\, \mathbb{R}^N \times (0,T), 
	\qquad 
	u(x,0) = u_0(x) 
	\quad\mbox{in} \,\,\, \mathbb{R}^N, 
\end{equation} 
satisfying $\displaystyle \lim_{t\to 0}\|u(t)-e^{t\Delta}u_0\|_{L^{\infty}(\mathbb{R}^N)}=0$.
\end{Proposition}
%
%

For the proof of Proposition~\ref{Proposition:2.5}, we introduce one lemma. 
\begin{Lemma} 
\label{Lemma:2.5}
Let $g$ be a convex function in $(s_0, \infty)$ 
such that 
$g(s)\ge Cs$ for all $s\ge s_0$ with some $C>0$ and $s_0>0$. 
Let $u_0\in L_{ul, \rho}^1(\mathbb{R}^N)$ be such that $u_0(x)\ge s_0$ in $\mathbb{R}^N$. 
Assume that there exists a nonnegative solution $u\in C^{2,1}(\mathbb{R}^N \times (0,T))$
satisfying 
$$
 u(x,t)=e^{t\Delta}u_0 +\int_0^t e^{(t-s)\Delta}e^{g(u(s))}ds \quad\mbox{in}\,\,\, \mathbb{R}^N\times (0,T).
$$
Then, for any $k\in \mathbb{N}$ with $k\ge 2$, there exists a constant $C_k>0$ such that
$$
 \|e^{t\Delta}u_0\|_{L^{\infty}(\mathbb{R}^N)}
 \le 
 g^{-1}\left(\frac{k}{k-1}
 \log \frac{1}{t}
 +C_k
 \right)
  \quad \mbox{for all}\,\,\, t\in (0,T).
$$
\end{Lemma}
\noindent 
{\bf Proof.} 
The proof of Lemma~\ref{Lemma:2.5} relies on the iteration argument developed by Weissler in \cite{W1}. 
Let $k\in \mathbb{N}$ and $k\ge 2$.  
We first prove by induction that 
\begin{equation} 
 \label{eq:2.24}
 u(x,t)\ge 
 \frac{
 \displaystyle
 t^{\displaystyle{a_l}}C^{\frac{k^l-1}{k-1}k}\exp\left({\displaystyle k^l\cdot g(e^{t\Delta}u_0)}\right)
 }
 {
 \displaystyle
 \prod_{i=1}^l (k!a_i)^{\displaystyle k^{l-i}}
 }
\end{equation}
for all $l\in \mathbb{N}$, $x\in \mathbb{R}^N$ and $t\in (0,T)$,
where $\{a_l\}$ is a sequence defined by 
$a_{l+1} = ka_l + 1$ with  $a_1 =  k + 1$. 

We start with the proof for the case $l=1$. 
Since $u$ satisfies the integral equation, we have $u(x,t)\ge e^{t\Delta}u_0(x)\ge s_0$ in $\mathbb{R}^N\times (0,T)$. 
Furthermore, since $u(x,t)\ge 0$ satisfies the integral equation and $\eta \in \mathbb{R}\mapsto e^{g(\eta)}\in \mathbb{R}$ is convex,  
by Lemma~\ref{Lemma:2.4} we have 
$$
 u(x,t)  \ge  \int_0^t e^{(t-s)\Delta}\exp\left(g(e^{s\Delta}u_0)\right) ds  
 \ge \int_0^t\exp\left(g(e^{(t-s)\Delta}e^{s\Delta}u_0)\right) ds=te^{g(e^{t\Delta}u_0)}.
$$
It follows from $e^{g(u)}\ge \frac{g(u)^{k}}{k!}\ge \frac{C^k u^k}{k!}$ for $u\ge s_0$, 
the convexity of $\eta\in \mathbb{R}\mapsto e^{kg(\eta)}\in \mathbb{R}$ and Lemma~\ref{Lemma:2.4} that
\begin{equation*} 
 \begin{split}
  u(x,t)
  &
  \ge 
  \int_0^t e^{(t-s)\Delta}\left(\frac{C^ku(s)^{k}}{k!}\right)\, ds
  \ge 
  \frac{C^k}{k!}\int_0^t e^{(t-s)\Delta}\left(s^{k}e^{kg(e^{s\Delta}u_0)}\right) ds
  \\
  &
  \ge 
  \frac{C^k}{k!}\int_0^t s^{k}e^{kg(e^{t\Delta}u_0)}\, ds 
  =
  \frac{C^k}{k!(k+1)}t^{k+1}e^{kg(e^{t\Delta}u_0)}.
 \end{split}
\end{equation*} 
This proves the inequality \eqref{eq:2.24} for the case $l=1$.

Now we assume that \eqref{eq:2.24} holds for $l\in \mathbb{N}$.
Applying $e^{g(u)}\ge \frac{C^ku^{k}}{k!}$ and \eqref{eq:2.24} with $l$, 
we have
\begin{equation*} 
  u(x,t)\ge 
  \frac{C^k}{k!}
  \int_0^t e^{(t-s)\Delta} 
  \left(
    \frac{
 \displaystyle
 s^{\displaystyle{ka_l}}C^{\frac{k^l-1}{k-1}k^2}\exp\left({\displaystyle k^{l+1}g(e^{s\Delta}u_0)}\right)
 }
 {
 \displaystyle
 \prod_{i=1}^l (k!a_i)^{\displaystyle k^{l-i+1}}
 }
  \right)
  ds.
\end{equation*} 
Again, 
by Lemma~\ref{Lemma:2.4} for the function $\eta \mapsto e^{k^{l+1}g(\eta)}$ 
we have 
\begin{align*} 
 u(x,t)
 &
 \ge 
  \frac{C^{\frac{k^{l+1}-1}{k-1}k}}{k!}
  \int_0^t 
    \frac{
 \displaystyle
 s^{\displaystyle{ka_l}}\exp\left({\displaystyle k^{l+1}g(e^{t\Delta}u_0)}\right)
 }
 {
 \displaystyle
 \prod_{i=1}^l (k!a_i)^{\displaystyle k^{l-i+1}}
 }
 \, ds 
 \\
 &
 =
    \frac{
 \displaystyle
 t^{\displaystyle{ka_l+1}} C^{\frac{k^{l+1}-1}{k-1}k}\exp\left({\displaystyle k^{l+1}g(e^{t\Delta}u_0)}\right)
 }
 {
 \displaystyle
 k!(ka_l+1) \prod_{i=1}^l (k!a_i)^{\displaystyle k^{l-i+1}}
 } 
 =
    \frac{
 \displaystyle
 t^{\displaystyle{a_{l+1}}} C^{\frac{k^{l+1}-1}{k-1}k}\exp\left({\displaystyle k^{l+1}g(e^{t\Delta}u_0)}\right)
 }
 {
 \displaystyle
 \prod_{i=1}^{l+1} (k!a_i)^{\displaystyle k^{l-i+1}}
 }. 
\end{align*} 
Here we used the relation $a_{l+1}=ka_l + 1$. 
This implies that \eqref{eq:2.24} holds with $l+1$. 
Thus we complete the proof of \eqref{eq:2.24}.

We now prove Lemma~\ref{Lemma:2.5}. 
It is easy to see that
$$
 a_l=\frac{k}{k-1}\cdot k^l -\frac{1}{k-1}
$$
for all $l\in \mathbb{N}$. 
Therefore, it follows from \eqref{eq:2.24} that
$$
 u(x,t)^{\frac{1}{k^l}}\prod_{i=1}^{l} (k!a_i)^{k^{-i}} 
 \ge 
 t^{\frac{k}{k-1}-\frac{1}{k^l(k-1)}}C^{\frac{k}{k-1}-\frac{k}{k^l(k-1)}}
 e^{g(e^{t\Delta}u_0)}
$$ 
for all $l\in \mathbb{N}$. Taking $l\to \infty$,
we have
\begin{equation} 
\label{eq:2.25}
 \prod_{i=1}^{\infty}(k!a_i)^{k^{-i}}\ge t^{\frac{k}{k-1}}C^{\frac{k}{k-1}}e^{g(e^{t\Delta}u_0)}.
\end{equation}
Remark that the left hand side of \eqref{eq:2.25} converges.
Indeed, we can easily see that 
$$
 \log\left(\displaystyle \prod_{i=1}^{\infty}(k!a_i)^{k^{-i}}\right)
 \le 
 \sum_{i=1}^{\infty} \frac{(i+1)\log k+\log (k!)}{k^i} <\infty. 
$$
Then we obtain the assertion of the lemma from \eqref{eq:2.25}, and complete the proof of Lemma~\ref{Lemma:2.5}.  
\hfill \qed 

\vspace{5pt} 

\noindent 
{\bf Proof of Proposition~\ref{Proposition:2.5}.} 
The proof is by contradiction. 
Fix $r<N/2$ and $2<\alpha <N/r$. 
Let $\epsilon>0$ be a sufficiently small constant such that $(1+\epsilon) \alpha r<N$. 
Since $g'>0$ and $g$ is convex, 
taking a sufficiently large $s_0>0$ if necessary, 
we can take a constant $C>0$ such that 
$g(s)\ge Cs$ for all $s\ge s_0$. 
Define
$$
u_0(x):=
\left\{
\begin{array}{lll} 
g^{-1}\left( 
\alpha \log\dfrac{1}{|x|} 
\right) 
&\mbox{if} &|x|<r_0,\\
s_0 &\mbox{if} &|x|\ge r_0, 
\end{array}
\right.
$$
where $r_0>0$ is chosen to satisfy $g^{-1}(\alpha\log (1/r_0)) = s_0$. 
Then, by $(1+\epsilon) \alpha r<N$ we have $e^{(1+\epsilon)r g(u_0)}\in L_{ul, \rho}^1(\mathbb{R}^N)$. 
Taking a sufficiently large $s_0>0$ if necessary, 
by \eqref{eq:2.22} we may assume that $g''(s)\le \epsilon (g'(s))^2$ for all $s\ge s_0$. 
Then we can easily see that 
\begin{equation} 
\label{eq:2.26} 
g'(s)\le g'(s_0) e^{\epsilon (g(s)-g(s_0))} 
\le g'(s_0) e^{\epsilon g(s)} 
\end{equation} 
for all $s\ge s_0$. 
On the other hand, since it follows from \eqref{eq:2.22} that 
$$
\frac{g''(s)}{(g'(s))^2}e^{-g(s)} \le \frac{1}{2} \cdot \frac{g''(s)}{(g'(s))^2} e^{-g(s)} 
+\frac{1}{2}e^{-g(s)} 
=\frac{1}{2}\left( 
-\frac{1}{g'(s)}e^{-g(s)} 
\right)' 
$$
for sufficiently large $s>0$, we have 
$$
\int_s^\infty \frac{g''(u)}{(g'(u))^2}e^{-g(u)}\, du 
\le \frac{1}{2}\cdot \frac{1}{g'(s)}e^{-g(s)}, 
$$
and by \eqref{eq:2.26} we obtain 
\begin{align*} 
G(s) 
& 
=\int_s^\infty \frac{du}{e^{g(u)}} 
=\frac{1}{g'(s)}e^{-g(s)} - \int_s^\infty \frac{g''(u)}{(g'(u))^2}e^{-g(u)}\, du 
\\ 
&
\ge 
\frac{1}{2}\cdot \frac{1}{g'(s)}e^{-g(s)} 
\ge 
\frac{1}{2}\cdot \frac{1}{g'(s_0)} e^{-(1+\epsilon) g(s)} 
\end{align*} 
for all sufficiently large $s>0$. 
This together with $e^{(1+\epsilon) g(u_0)}\in L_{ul, \rho}^1(\mathbb{R}^N)$ implies that $G(u_0)^{-r}\in L^1_{ul, \rho}(\mathbb{R}^N)$. 

Putting $y=\sqrt{t}z$, we have 
\begin{align*} 
	\|e^{t\Delta}u_0\|_{L^\infty(\mathbb{R}^N)} 
	& 
	\ge (4\pi t)^{-\frac{N}{2}} 
	\int_{|y|\le r_0} e^{-\frac{|y|^2}{4t}} g^{-1}\left( \alpha \log \frac{1}{|y|}\right) dy 
	\\ 
	&
	= (4\pi)^{-\frac{N}{2}} 
	\int_{|z|\le r_0 t^{-1/2}} e^{-\frac{|z|^2}{4}} 
	g^{-1} \left( 
	\frac{\alpha}{2}\log \frac{1}{t} + \alpha\log \frac{1}{|z|} 
	\right) dz 
\end{align*} 
Let $\delta>0$ be a sufficiently small constant such that $\alpha/2-\delta>1+2\delta$. 
Then we have $t^{-\delta/\alpha}\le r_0 t^{-1/2}$ for all sufficiently small $t>0$ and 
\begin{equation} 
\label{eq:2.27} 
\begin{split} 
	\|e^{t\Delta}u_0\|_{L^\infty(\mathbb{R}^N)} 
	&
	\ge (4\pi)^{-\frac{N}{2}} 
	\int_{|z|\le t^{-\delta/\alpha}} e^{-\frac{|z|^2}{4}} 
	g^{-1} \left( 
	\frac{\alpha}{2}\log \frac{1}{t} + \alpha\log \frac{1}{|z|} 
	\right) dz 
	\\ 
	&
	\ge g^{-1}\left( 
	\left(\frac{\alpha}{2}-\delta\right) \log \frac{1}{t} 
	\right) 
	\cdot 
	(4\pi)^{-\frac{N}{2}}\int_{|z|\le t^{-\delta/\alpha}} e^{-\frac{|z|^2}{4}}\, dz 
	\\ 
	&
	=  g^{-1}\left( 
	\left(\frac{\alpha}{2}-\delta\right) \log \frac{1}{t} 
	\right) 
	\cdot 
	\left( 1+O(e^{-\frac{t^{-2\delta/\alpha}}{8}}) \right) 
\end{split} 
\end{equation} 
for all sufficiently small $t>0$. 
We prove that this yields a contradiction. 
Assume that there exists a nonnegative classical solution $u$ of 
$$
	u(t) = e^{t\Delta} u_0 + \int_0^t e^{(t-s)\Delta} e^{g(u(s))}\, ds. 
$$
Then, by Lemma~\ref{Lemma:2.5} and \eqref{eq:2.27} we have  
\begin{equation} 
\label{eq:2.28} 
	g^{-1}\left( 
	\left(\frac{\alpha}{2}-\delta\right) \log \frac{1}{t} 
	\right) 
	\cdot 
	\left( 1+ O(e^{-\frac{t^{-2\delta/\alpha}}{8}}) \right) 
	\le 
	g^{-1}\left( 
	\left(\frac{k}{k-1} + \delta\right)\log \frac{1}{t} 
	\right) 
\end{equation} 
for all sufficiently small $t>0$. 
Here we take a sufficiently large $k\in \mathbb{N}$ satisfying $k/(k-1)+\delta<1+2\delta$. 
Then, since $\alpha/2-\delta>1+2\delta$, we have $k/(k-1)+\delta < \alpha/2-\delta$. 
In the following, we prove that \eqref{eq:2.28} yields a contradiction. 
For simplicity, define $a:=\alpha/2-\delta$, $b:=k/(k-1)+\delta$, $c:=2\delta/\alpha$ and $\tau := \log (1/t)$. 
Then $a>b$ and \eqref{eq:2.28} implies that 
\begin{equation} 
\label{eq:2.29} 
	g^{-1}(a\tau) \cdot 
	\left(1+O(e^{-e^{-c\tau}/8}) \right) 
	\le 
	g^{-1}(b\tau) 
\end{equation} 
for sufficiently large $\tau>0$. 
Since $(g^{-1}(s))'=1/g'(g^{-1}(s))>0$,  
by the mean value theorem and \eqref{eq:2.26} we see that there exists a constant $d>0$ such that 
$$
	g^{-1}(a\tau) 
	\ge g^{-1}(b\tau)+\frac{1}{g'(g^{-1}(a\tau))}(a-b)\tau 
	\ge g^{-1}(b\tau)+\frac{d}{e^{\epsilon a\tau}}(a-b)\tau, 
$$
which implies that 
$$
	g^{-1}(a\tau) \cdot 
	\left(1+O(e^{-e^{-c\tau}/8}) \right) 
	\ge 
	g^{-1}(b\tau)+\frac{d}{e^{\epsilon a \tau}}(a-b)\tau + O(g^{-1}(a\tau) e^{-e^{-c\tau}/8}) 
	>g^{-1}(b\tau) 
$$
for all sufficiently large $\tau>0$. 
Remark that $g^{-1}(a\tau)=O(\tau)$ for all sufficiently large $\tau$ since $g(s) \ge Cs$ for all $s\ge s_0$. 
This contradicts \eqref{eq:2.29}, and so \eqref{eq:2.28} yields a contradiction. 

We finally prove the latter assertion of Proposition~\ref{Proposition:2.5}. 
Assume that there exist a constant $T>0$ and 
a nonnegative classical solution $u\in C^{2,1}(\mathbb{R}^N \times (0,T))$ 
of \eqref{eq:2.23} satisfying $\displaystyle \lim_{t\to 0}\|u(t)-e^{t\Delta}u_0\|_{L^{\infty}(\mathbb{R}^N)}=0$. 
Then, for any $\tau\in (0,T)$, by \eqref{eq:2.23} we have 
\begin{equation} 
\label{eq:2.30} 
u(x,t) = (e^{(t-\tau)\Delta}u(\tau))(x)+\int_\tau^t [e^{(t-s)\Delta}e^{g(u(s))}](x)\, ds 
\end{equation} 
in $\mathbb{R}^N\times (\tau, T)$. 
This implies that 
$$
\int_\tau^t [e^{(t-s)\Delta}e^{g(u(s))}](x)\, ds \le u(x,t)<\infty, 
$$
so we see that 
$$
\lim_{\tau\to 0}\int_\tau^t [e^{(t-s)\Delta}e^{g(u(s))}](x)\, ds =\int_0^t [e^{(t-s)\Delta}e^{g(u(s))}](x)\, ds 
$$
exists for all $(x,t)\in\mathbb{R}^N\times (0,T)$. 
Therefore, since $\|u(\tau)-e^{\tau\Delta}u_0\|_{L^{\infty}(\mathbb{R}^N)}\to 0$ as $\tau\to 0$,  
as in the proof of Lemma~\ref{Lemma:2.3}, 
taking the limit $\tau\to 0$ in \eqref{eq:2.30}, 
we see that $u$ satisfies the integral equation
$$
 u(x,t)=e^{t\Delta}u_0 +\int_0^t e^{(t-s)\Delta}e^{g(u(s))}ds \quad\mbox{in}\,\,\, \mathbb{R}^N  \times (0,T).
$$
This is a contradiction. 
Thus we complete the proof of Proposition~\ref{Proposition:2.5}. 
\hfill 
\qed 

\section{Existence of a solution for problem~\eqref{eq:1.1}}
\label{section:3} 
In this section we show local and global in time existence of solutions for problem~\eqref{eq:1.1} 
with the help of Propositions~\ref{Proposition:2.1}--\ref{Proposition:2.4}. 
Recall
that
\begin{equation} 
\label{eq:3.1} 
A=\lim_{s\to\infty}f'(s)F(s). 
\end{equation} 
Before starting the proof of main theorems, 
we prepare two lemmas.
%
\begin{Lemma} 
\label{Lemma:3.1} 
Let $f\in C^1([0,\infty))$ satisfy \eqref{eq:1.2} and \eqref{eq:1.5}. 
Define $g(s):= f\left(F^{-1}(s)\right)$. 
Assume that \eqref{eq:1.9} holds for some $s_1>0$. 
Then 
\begin{equation} 
\notag 
 g(s)\lesssim s^{-A} 
\end{equation}
for all sufficiently small $s>0$. 
\end{Lemma}
\noindent 
{\bf Proof.} 
It follows from $(F^{-1}(s))'=-f(F^{-1}(s))=-g(s)$ that
\begin{equation}\label{eq:3.2}
 g'(s)=f'(F^{-1}(s))(F^{-1}(s))'=-f'(F^{-1}(s))g(s).
\end{equation}
Since $F^{-1}$ is monotonically decreasing with respect to $s$ and $F^{-1}(s)\to \infty$ as $s\to 0$, 
we have $F^{-1}(s)>s_1$ for all sufficiently small $s>0$. 
This together with \eqref{eq:1.9} implies that 
\begin{equation}\label{eq:3.3}
 sf'(F^{-1}(s))=f'(F^{-1}(s))F(F^{-1}(s))\le A 
\end{equation}
for all sufficiently small $s>0$. 
Combining \eqref{eq:3.2} and \eqref{eq:3.3}, we have
$$
g'(s)\ge -As^{-1}g(s) 
$$
for all sufficiently small $s>0$. 
This proves Lemma~\ref{Lemma:3.1}. 
\hfill \qed 

\begin{Lemma}　
\label{Lemma:3.2} 
Let $f\in C^1([0,\infty))$ satisfy \eqref{eq:1.2} and \eqref{eq:1.5}. 
Let $A>1$ and $s_1>0$. 
\begin{itemize}
\item[\rm (i)]
If $f'(s)F(s)\le A$ for all $s>s_1$, then
$F(s)^{-(A-1)}$ is convex in $(s_1,\infty)$.
On the other hand, 
if
$f'(s)F(s)\ge A$ for all $s>s_1$,
then $F(s)^{-(A-1)}$ is concave in $(s_1,\infty)$.
\item[\rm (ii)]
If
$f'(s)F(s)\le 1$ for all $s>s_1$,
then
$\log F(s)^{-1}$ is convex in $(s_1,\infty)$.
On the other hand, 
if
$f'(s)F(s)\ge 1$ for all $s>s_1$,
then $\log F(s)^{-1}$ is concave in $(s_1,\infty)$.
\end{itemize}
\end{Lemma} 
\noindent 
{\bf Proof.} 
Assume that $f'(s)F(s)\le A$ for all $s>s_1$. 
Then direct calculations show that
$$
 \frac{d^2}{ds^2}\left(F(s)^{-(A-1)}\right)=(A-1)F(s)^{-A-1}\frac{1}{f(s)^2}\left(A-f'(s)F(s)\right)\ge 0 
$$
for all $s>s_1$. 
This proves the convexity of $F(s)^{-(A-1)}$ in assertion~(i).
Other cases can be treated in the same manner.
\hfill \qed 

\vspace{5pt} 

%
%
The rest of this section is devoted to the proof of Theorem~\ref{Theorem:1.1}. 
\vspace{5pt} 
\newline 
{\bf Proof of Theorem~\ref{Theorem:1.1}.} 
\vspace{5pt} 
\newline 
\underline{\bf Case $A>1$} 
\vspace{3pt} 
\newline 
We first consider the case  $A>1$, and prove local in time existence of a solution for \eqref{eq:1.1}.
Let $r$ be a constant given in the assumption of Theorem~\ref{Theorem:1.1}. 
Define
\begin{equation} 
\label{eq:3.4} 
v_0(x):=\max\left\{ 
F(u_0(x))^{-(A-1)}, F(s_1)^{-(A-1)} 
\right\},
\end{equation} 
where $s_1$ is the constant appearing in \eqref{eq:1.9}.
In particular, we have 
\begin{equation} 
\label{eq:3.5} 
v_0(x)\ge F(s_1)^{-(A-1)}>0 \quad\mbox{in}\,\,\, {\mathbb{R}}^N. 
\end{equation} 
Consider the semilinear heat equation 
\begin{equation} 
\label{eq:3.6} 
\partial_t v=\Delta v+(A-1)v^{\frac{A}{A-1}} \quad \mbox{in} \,\,\, \mathbb{R}^N\times (0,T), \qquad 
v(x,0)=v_0(x) \quad \mbox{in}\,\,\, \mathbb{R}^N, 
\end{equation} 
where $T>0$. 
By \eqref{eq:3.4} we have 
\begin{equation} 
\label{eq:3.7}
\sup_{y\in {\mathbb{R}}^N}\int_{B_{\rho}(y)} |v_0(x)|^{\frac{r}{A-1}}\, dx 
=  
\sup_{y\in {\mathbb{R}}^N}\int_{B_{\rho}(y)} 
\max\left\{F(u_0(x))^{-r}, F(s_1)^{-r}\right\} dx, 
\end{equation} 
which implies that $v_0\in {L}_{ul,\rho}^{r/(A-1)}({\mathbb{R}}^N)$ in case of \eqref{eq:1.11}  
and 
$v_0\in {\cal L}_{ul,\rho}^{r/(A-1)}(\mathbb{R}^N)$ in case of \eqref{eq:1.17}. 
Since 
\begin{itemize} 
\item 
$\dfrac{r}{A-1}>\dfrac{N}{2(A-1)}=\dfrac{N}{2}\left[\dfrac{A}{A-1}-1\right]$ and $\dfrac{r}{A-1}\ge 1$ 
if \eqref{eq:1.10} holds; 
\item 
$\dfrac{r}{A-1}\ge \dfrac{N}{2(A-1)}=\dfrac{N}{2}\left[\dfrac{A}{A-1}-1\right]$ and $\dfrac{r}{A-1}>1$ 
if \eqref{eq:1.16} holds,  
\end{itemize} 
in view of Propositions~\ref{Proposition:2.2} and \ref{Proposition:2.3}, 
there exist a constant $T>0$ and a classical solution $v$ in $L^{r/(A-1)}_{ul,\rho}(\mathbb{R}^N)$ (resp. in ${\cal L}_{ul,\rho}^{r/(A-1)}({\mathbb{R}}^N)$) 
of \eqref{eq:3.6} satisfying 
\begin{equation} 
\label{eq:3.8} 
	v(x,t) = (e^{t\Delta}v_0)(x) 
	+ \int_0^t [e^{(t-s)\Delta}v(s)^\frac{A}{A-1}](x) \, ds 
	\quad\mbox{in}\,\,\, 
	\mathbb{R}^N \times (0,T), 
\end{equation} 
if $r$ satisfies \eqref{eq:1.10} (resp. $r$ satisfies \eqref{eq:1.16}). 
Then we have $v\in C^{2,1}({\mathbb{R}}^N\times (0,T))$. 
Furthermore, by \eqref{eq:3.5} and \eqref{eq:3.8} we obtain 
\begin{equation} 
\label{eq:3.9} 
v(x,t)
\ge (e^{t\Delta}v_0)(x) \ge F(s_1)^{-(A-1)} 
\end{equation} 
in $\mathbb{R}^N\times (0,T)$. 

Define the function $\overline{u}\in C^{2,1}({\mathbb{R}}^N\times (0,T))$ by 
\begin{equation} 
\label{eq:3.10}
\overline{u}(x,t):=F^{-1}\left(v(x,t)^{-1/(A-1)}\right). 
\end{equation}
Then, by \eqref{eq:3.9} we have 
\begin{equation} 
\label{eq:3.11} 
\overline{u}(x,t)\ge F^{-1}(F(s_1))= s_1 \quad\mbox{in}\,\,\, {\mathbb{R}}^N\times (0,T). 
\end{equation} 
Since we have 
\begin{equation} 
\label{eq:3.12} 
\begin{array}{l} 
\partial_t\overline{u}=\dfrac{1}{A-1}f(\overline{u})v^{-\frac{A}{A-1}}\partial_t v, \vspace{3pt} \\ 
\Delta\overline{u}=\dfrac{f'(\overline{u})f(\overline{u})v^{-\frac{2A}{A-1}}|\nabla v|^2}{(A-1)^2} 
-\dfrac{Af(\overline{u})v^{-\frac{A}{A-1}-1}|\nabla v|^2}{(A-1)^2}  
+\dfrac{f(\overline{u})v^{-\frac{A}{A-1}}\Delta v}{A-1}, 
\end{array} 
\end{equation} 
by \eqref{eq:3.6} we obtain 
\begin{align*} 
\partial_t\overline{u}-\Delta \overline{u}-f(\overline{u}) 
&=\frac{1}{(A-1)^2}f(\overline{u})v^{-\frac{A}{A-1}-1}|\nabla v|^2\Big[ 
A-f'(\overline{u})v^{-\frac{1}{A-1}}\Big]
\\ 
&=\frac{1}{(A-1)^2}f(\overline{u})v^{-\frac{A}{A-1}-1}|\nabla v|^2\Big[ 
A-f'(\overline{u})F(\overline{u})\Big].
\end{align*} 
Thus, since $f'(\overline{u})F(\overline{u})\le A$ in $\mathbb{R}^N\times (0,T)$ by \eqref{eq:1.9} and \eqref{eq:3.11}, we obtain 
\begin{equation} 
\label{eq:3.13} 
\partial_t\overline{u}\ge \Delta\overline{u}+f(\overline{u}) \quad\mbox{in}\,\,\, {\mathbb{R}}^N\times (0,T). 
\end{equation} 
Furthermore, by \eqref{eq:3.4} we see that  
\begin{equation} 
\label{eq:3.14} 
\overline{u}(x,0)=\max\left\{u_0(x), s_1\right\} \ge u_0(x) \quad\mbox{in}\,\,\, {\mathbb{R}}^N. 
\end{equation} 
By \eqref{eq:3.10} and \eqref{eq:3.11} we obtain $F^{-1}(v(x,t)^{-1/(A-1)})\ge s_1$. 
This together with Lemma~\ref{Lemma:3.1} with $A>1$ implies that 
\begin{equation} 
\label{eq:3.15} 
  \int_0^t e^{(t-s)\Delta}f(\overline{u}(s)) \, ds 
  = \int_0^t e^{(t-s)\Delta}f(F^{-1}(v(s)^{-\frac{1}{A-1}})) \, ds 
  \lesssim \int_0^t e^{(t-s)\Delta}v(s)^{\frac{A}{A-1}} \, ds 
\end{equation} 
for all $(x,t)\in \mathbb{R}^N\times (0,T)$. 
For the case $r=N/2$, since $\frac{r}{A-1}=\frac{N}{2}(\frac{A}{A-1}-1)$ and 
$v\in C([0,T]; {\cal L}_{ul,\rho}^{r/(A-1)}(\mathbb{R}^N))$, 
we have 
\begin{equation} 
\label{eq:3.16} 
  \int_0^t e^{(t-s)\Delta}v(s)^{\frac{A}{A-1}} \, ds 
=
v(t)-e^{t\Delta}v_0  \to 0 \quad\mbox{in} \,\,\, L^{\frac{r}{A-1}}_{ul,\rho}(\mathbb{R}^N), 
\end{equation} 
as $t\to 0$. 
For the case $r>N/2$, by Proposition~\ref{Proposition:2.2} we also have \eqref{eq:3.16}. 
Then, since $\frac{r}{A-1}\ge 1$, by \eqref{eq:3.15} and \eqref{eq:3.16} we have \eqref{eq:2.3} for $\overline{u}$. 
On the other hand, 
by Lemma~\ref{Lemma:2.4}, Lemma~\ref{Lemma:3.2}, \eqref{eq:1.9} and \eqref{eq:3.14} we have 
$$
F(e^{t\Delta}\overline{u}(0))^{-(A-1)}\le e^{t\Delta}\left[F(\overline{u}(0))^{-(A-1)}\right]. 
$$
This together with \eqref{eq:3.9} and \eqref{eq:3.10} implies that 
$$
F(\overline{u}(t))^{-(A-1)}-F(e^{t\Delta}\overline{u}(0))^{-(A-1)} 
\ge v(t)-e^{t\Delta}\left[ F(\overline{u}(0))^{-(A-1)} \right] 
=v(t)-e^{t\Delta}v_0\ge 0. 
$$
Since $F(s)^{-(A-1)}$ is monotonically increasing with respect to $s$, we see that $\overline{u}(x,t)\ge (e^{t\Delta}\overline{u}(0))(x)$ in 
$\mathbb{R}^N\times (0,T)$. 
Therefore we can apply Lemma~\ref{Lemma:2.3}~(ii)
and by \eqref{eq:3.13} and \eqref{eq:3.14}
we obtain 
$$
\overline{u}(x,t)\ge (e^{t\Delta}\overline{u}(0))(x)+\int_0^t [e^{(t-s)\Delta}f(\overline{u}(s))](x)\, ds 
\ge (e^{t\Delta}u_0)(x)+\int_0^t [e^{(t-s)\Delta}f(\overline{u}(s))](x)\, ds 
$$
for all $(x,t)\in {\mathbb{R}}^N\times (0,T)$. 
Hence we can apply Proposition~\ref{Proposition:2.1} and prove existence of a local in time classical solution $u$ for problem~\eqref{eq:1.1}. 
Note that, in case of \eqref{eq:1.10}, the existence time $T$ satisfies 
\begin{equation}
\notag 
\begin{split} 
T^{\frac{N}{2}(1-\frac{1}{A})}\rho^{-N(1-\frac{1}{A})} 
+ 
\gamma_0 \max\left\{\|F(u_0)^{-r}\|_{L_{ul,\rho}^1(\mathbb{R}^N)}, F(s_1)^{-r}\rho^N \right\} 
\left( 
T^{r-\frac{N}{2A}}\rho^{-N\cdot\frac{A-1}{A}}+T^{r-\frac{N}{2}} 
\right) \\ 
\ge 
T^{\frac{N}{2}(1-\frac{1}{A})}\rho^{-N(1-\frac{1}{A})} 
+ \|v_0\|_{L^{\frac{r}{A-1}}_{ul,\rho}(\mathbb{R}^N)}^{\frac{r}{A-1}} 
\left( 
T^{r-\frac{N}{2A}}\rho^{-N\cdot\frac{A-1}{A}}+T^{r-\frac{N}{2}} 
\right) 
\ge 
\gamma_1 
\end{split} 
\end{equation}
by \eqref{eq:2.11} with $p=\frac{A}{A-1}$, \eqref{eq:3.4} and \eqref{eq:3.7}, 
where $\gamma_0$ and $\gamma_1$ are positive constants depending only on $N$, $A$ and $r$. 
Thus we obtain \eqref{eq:1.13}. 

It remains to prove \eqref{eq:1.12}.
By Proposition~\ref{Proposition:2.1} we have $u(x,t)\le \overline{u}(x,t)$. 
Then, by \eqref{eq:1.2} and \eqref{eq:3.15} we have 
\begin{equation} 
\notag 
  |u(t)-e^{t\Delta}u_0|
  =\int_0^t e^{(t-s)\Delta}f(u(s))ds\le  \int_0^t e^{(t-s)\Delta}f(\overline{u}(s)) \, ds  
  \lesssim \int_0^t e^{(t-s)\Delta}v(s)^{\frac{A}{A-1}} \, ds. 
\end{equation} 
This together with \eqref{eq:3.16} yields \eqref{eq:1.12}. 
Thus we complete the proof of Theorem~\ref{Theorem:1.1} for the case $A>1$. 
\vspace{5pt} 
\newline 
\underline{\bf Case $A=1$} 
\vspace{3pt} 
\newline 
We next consider the case $A=1$. 
The proof is similar to the above argument. 
Let 
\begin{equation} 
\label{eq:3.17} 
w_0(x):=\max\left\{ 
{\log F(u_0(x))^{-1}}, {\log F(s_1)^{-1}} 
\right\}.
\end{equation} 
Then we have 
\begin{equation} 
\label{eq:3.18} 
w_0(x)\ge {\log F(s_1)^{-1}} \quad\mbox{in}\,\,\, {\mathbb{R}}^N. 
\end{equation} 
Consider the semilinear heat equation 
\begin{equation} 
\label{eq:3.19} 
\partial_t w=\Delta w+e^w \quad \mbox{in}\,\,\, {\mathbb{R}}^N\times (0,T), \qquad 
w(x,0)=w_0(x) \quad \mbox{in} \,\,\, {\mathbb{R}}^N, 
\end{equation} 
where $T>0$. 
Since \eqref{eq:1.11} and \eqref{eq:3.17} yields 
\begin{equation} \label{eq:3.20}
\sup_{y\in \mathbb{R}^N}\int_{B_{\rho}(y)} e^{rw_0(x)}\, dx 
=
\sup_{y\in {\mathbb{R}}^N}\int_{B_{\rho}(y)}\max\left\{F(u_0(x))^{-r}, F(s_1)^{-r}\right\} dx, 
\end{equation} 
we have $e^{rw_0}\in L^1_{ul,\rho}(\mathbb{R}^N)$ if $r>N/2$ and $e^{rw_0}\in {\cal L}^1_{ul,\rho}(\mathbb{R}^N)$ if $r=N/2$. 
Hence, by Proposition~\ref{Proposition:2.4} we can find a constant $T>0$ 
and a classical solution $w$ of \eqref{eq:3.19} 
satisfying 
\begin{equation} 
\notag 
w(t)=e^{t\Delta}w_0+\int_0^t e^{(t-s)\Delta}e^{w(s)}\, ds\ge e^{t\Delta}w_0 
\end{equation} 
and 
\begin{equation} 
\label{eq:3.21} 
\lim_{t\to 0}\|w(t)-e^{t\Delta}w_0\|_{L^\infty(\mathbb{R}^N)}
=\lim_{t\to 0}\left\| 
\int_0^t e^{(t-s)\Delta}e^{w(s)}\, ds 
\right\|_{L^\infty(\mathbb{R}^N)} =0. 
\end{equation} 

Define the function $\overline{u}\in C^{2,1}({\mathbb{R}}^N\times (0,T))$ by 
\begin{equation}\label{eq:3.22}
\overline{u}(x,t):=F^{-1}\left(e^{-w(x,t)}\right). 
\end{equation}
Then, as in \eqref{eq:3.11}, by \eqref{eq:3.18} we have 
\begin{equation} 
\label{eq:3.23} 
\overline{u}(x,t)\ge s_1 \quad\mbox{in}\,\,\, {\mathbb{R}}^N\times (0,T). 
\end{equation} 
One can easily check that
$$
\partial_t \overline{u}
=f(\overline{u})e^{-w}\partial_t w, \quad 
\Delta \overline{u}
=f'(\overline{u})f(\overline{u})e^{-2w}|\nabla w|^2-f(\overline{u})e^{-w}|\nabla w|^2 +f(\overline{u})e^{-w}\Delta w, 
$$
and so we have  
\begin{equation} 
\notag 
\partial_t \overline{u}-\Delta \overline{u}-f(\overline{u})
=f(\overline{u})e^{-w}|\nabla w|^2( 1-f'(\overline{u})e^{-w})
=f(\overline{u})e^{-w}|\nabla w|^2( 1-f'(\overline{u})F(\overline{u})). 
\end{equation} 
Since $f'(\overline{u})F(\overline{u})\le 1$ by \eqref{eq:1.9} with $A=1$ and \eqref{eq:3.23}, we have
\begin{equation} 
\label{eq:3.24} 
\partial_t \overline{u}\ge \Delta \overline{u}+f(\overline{u}) \quad\mbox{in}\,\,\, {\mathbb{R}}^N\times (0,T). 
\end{equation} 
By the similar argument as in the case $A>1$ with the aid of Lemma~\ref{Lemma:3.2} we obtain $\overline{u}(x,t)\ge (e^{t\Delta}\overline{u}(0))(x)$. 
Furthermore, by Lemma~\ref{Lemma:3.1} with $A=1$ we have 
$$
\int_0^t e^{(t-s)\Delta}f(\overline{u}(s)) \, ds
= \int_0^t e^{(t-s)\Delta}f(F^{-1}(e^{-w}))\, ds 
\lesssim \int_0^t e^{(t-s)\Delta}e^w\, ds. 
$$
Then, by 
\eqref{eq:3.21} 
and 
\eqref{eq:3.24}, 
we can apply Lemma~\ref{Lemma:2.3}~(ii), and by \eqref{eq:3.17} we obtain 
$$
\overline{u}(x,t)\ge e^{t\Delta}u_0+\int_0^t e^{(t-s)\Delta}f(\overline{u}(s))\, ds. 
$$
Applying Proposition~\ref{Proposition:2.1}, we obtain a local in time classical solution $u$ of problem~\eqref{eq:1.1}. 
Note that, if $r>\frac{N}{2}$,
by \eqref{eq:2.14}, \eqref{eq:3.17} and \eqref{eq:3.20} we see that the existence time $T$ can be taken to satisfy 
\begin{align*} 
 T^{\frac{\epsilon N}{2}}\rho^{-\epsilon N} 
 + C\gamma_\epsilon \max\left\{\|e^{ru_0}\|_{L_{ul,\rho}^1(\mathbb{R}^N)}, F(s_1)^{-r}\rho^N\right\} 
 \left( 
 T^{r-\frac{N}{2}(1-\epsilon)}\rho^{-\epsilon N}+T^{r-\frac{N}{2}} 
 \right) \quad \\ 
 \ge 
 T^{\frac{\epsilon N}{2}}\rho^{-\epsilon N} 
 + \gamma_\epsilon \|e^{rw_0}\|_{L_{ul,\rho}^1(\mathbb{R}^N)} 
 \left( 
 T^{r-\frac{N}{2}(1-\epsilon)}\rho^{-\epsilon N}+T^{r-\frac{N}{2}} 
 \right) 
 \ge \gamma_* 
\end{align*} 
for all sufficiently small $\epsilon>0$, 
where $C>0$ is a constant and $\gamma_*$ depends only on $N$ and $r$, 
and $\gamma_\epsilon$ is a constant depending on $N$, $r$ and $\epsilon$ 
such that $\gamma_\epsilon\to \infty$ as $\epsilon\to 0$. 
Thus we obtain \eqref{eq:1.15}. 

Convergence of the solution $u$  
to the initial data \eqref{eq:1.14} is proved as follows. 
Since $u(x,t)\le \overline{u}(x,t) =F^{-1}(e^{-w})$ and $F^{-1}(e^{-w}) \ge s_1$ by \eqref{eq:3.22} and \eqref{eq:3.23}, 
we can apply Lemma~\ref{Lemma:3.1} with $A=1$ and obtain 
\begin{align*} 
|u(t)-e^{t\Delta}u_0|
= \int_0^t e^{(t-s)\Delta}f({u}(s))\, ds \le \int_0^t e^{(t-s)\Delta}f(\overline{u}(s)) \, ds 
\lesssim \int_0^t e^{(t-s)\Delta}e^w\, ds. 
\end{align*} 
Then \eqref{eq:3.21} shows us the desired convergence. 
Thus we complete the proof of Theorem~\ref{Theorem:1.1} for the case $A=1$. 
\hfill \qed 
\begin{Remark} 
\label{Remark:3.1} 
%
We explain the structure of the transformations \eqref{eq:3.10} and \eqref{eq:3.22} 
used in the proof of Theorem~{\rm\ref{Theorem:1.1}} for the cases $A>1$ and $A=1$, respectively.
Let $f$ and $g$ satisfy \eqref{eq:1.2} and \eqref{eq:1.5}, and define $G$ by
$$
G(v):=\int_v^{\infty}\frac{ds}{g(s)}.
$$
Assume that $v$ satisfies
$\partial_t v -\Delta v=g(v)$, and consider the following general transformation 
\begin{equation} 
\label{eq:3.25} 
	\tilde u(x,t):=F^{-1}\left(G(v(x,t))\right).
\end{equation} 
By a simple calculation we see that $\tilde u$ satisfies
$$
   \partial_t \tilde{u}-  \Delta\tilde{u}-f(\tilde u)
   =\frac{f(\tilde{u})|\nabla v|^2}{g(v)^2F(\tilde u)}\left(g'(v)G(v)-f'(\tilde u)F(\tilde u)\right).
$$
Then one can easily check that 
$$
 g(s)=
\left\{
\begin{array}{lll} 
s^{\frac{A}{A-1}} &\text{if}  &A>1,
\vspace{3pt} \\
e^s &\text{if} &A=1,
\end{array}
\right.
$$
is a solution of the equation $g'(s)G(s)=A$, which implies that 
$G(s) = (A-1)s^{-\frac{1}{A-1}}$ for the case $A>1$ and $G(s) = e^{-s}$ for the case $A=1$. 
Then \eqref{eq:3.25} corresponds to \eqref{eq:3.10} for the case $A>1$ and 
\eqref{eq:3.22} for the case $A=1$, respectively. 
\end{Remark} 
%

\section{Nonexistence of solutions for problem~\eqref{eq:1.1}}
\label{section:4} 
In this section we discuss the nonexistence results of local in time nonnegative classical  
solutions of \eqref{eq:1.1}, 
and prove Theorems~\ref{Theorem:1.3} and \ref{Theorem:1.4}. 
Recall 
$$
A=\lim_{s\to \infty}f'(s)F(s). 
$$
{\bf Proof of Theorem~\ref{Theorem:1.3}.} 
We first consider the case $A>1$. 
Let $r\in [A-1, N/2)$. 
By \cite[Corollary~5.1]{W1} we can take a nonnegative function 
$v_0 \in L^{r/(A-1)}(\mathbb{R}^N)\subset {L}_{ul,\rho}^{r/(A-1)}(\mathbb{R}^N)$ such that
there can not exists a local in time solution for the integral equation 
\begin{equation} 
\label{eq:4.1} 
v(t)=e^{t\Delta}v_0+(A-1)\int_0^t e^{(t-s)\Delta}v(s)^\frac{A}{A-1}\, ds.
\end{equation} 
Remark that, under the assumption of Theorem~\ref{Theorem:1.3} for the case $A>1$, we have 
$$
1\le \frac{r}{A-1}<\frac{N}{2}\left[\frac{A}{A-1}-1\right]. 
$$
Let $s_2>0$ be the constant satisfying \eqref{eq:1.20} for all $s\ge s_2$. 
Define
\begin{equation}\label{eq:4.2}
u_0(x):=\max\left\{F^{-1}\left(v_0^{-\frac{1}{A-1}}\right),s_2\right\}.
\end{equation}
The proof is by contradiction. 
Suppose that there exist an existence time $T>0$ and  
a local in time nonnegative classical solution $\overline{u}\in C^{2,1}(\mathbb{R}^N\times (0,T))$ 
for the problem~\eqref{eq:1.1} satisfying
\begin{equation}\label{eq:4.3}
\lim_{t\to 0}\|u(t)-e^{t\Delta}u_0\|_{L_{ul,\rho}^{\frac{r}{A-1}}(\mathbb{R}^N)} =0.
\end{equation}
Note that $F(u_0)^{-r}\in {L}_{ul,\rho}^{1}(\mathbb{R}^N)$. 
Since $\overline{u}$ is a classical solution of \eqref{eq:1.1} satisfying \eqref{eq:4.3}, 
as in the proof of Proposition~\ref{Proposition:2.5}, 
we can rewrite \eqref{eq:1.1} by the integral form 
\begin{equation} 
\label{eq:4.4} 
\overline u(x,t)=e^{t\Delta} u_0+\int_0^t e^{(t-s)\Delta}f(\overline u(s))\, ds.
\end{equation} 
This implies $\overline u(x,t)\ge e^{t\Delta}u_0\ge s_2$
and
\begin{equation} 
\label{eq:4.5} 
f'(\overline u(x,t))F(\overline u(x,t))\ge A 
\end{equation} 
in $\mathbb{R}^N\times (0,T)$. 
Define 
\begin{equation} 
\notag 
 \overline v(x,t):=F(\overline u(x,t))^{-(A-1)} \quad\mbox{in}\,\,\, \mathbb{R}^N\times (0,T). 
\end{equation}
Then, by \eqref{eq:4.5} we see that $\overline{v}$ satisfies 
\begin{equation} 
\label{eq:4.6} 
\partial_t \overline v-\Delta \overline v-(A-1)\overline v^{\frac{A}{A-1}} 
=(A-1)F(\overline u)^{-A-1}f(\overline u)^{-2}|\nabla \overline u|^2\Big[f'(\overline u)F(\overline u)-A\Big]\ge 0 
\end{equation} 
in $\mathbb{R}^N\times (0,T)$. 
By Lemma~\ref{Lemma:3.2}~(i) together with \eqref{eq:1.20} we can apply Lemma~\ref{Lemma:2.4}, 
and by \eqref{eq:4.2} and \eqref{eq:4.4} we have 
\begin{equation} 
\label{eq:4.7} 
\overline{v}(x,t)\ge F(e^{t\Delta}u_0)^{-(A-1)}\ge e^{t\Delta}\left[F(u_0)^{-(A-1)}\right] 
\ge 
e^{t\Delta}v_0. 
\end{equation} 
Furthermore, for any $\tau\in (0,T)$, by \eqref{eq:4.6} we have 
\begin{equation} 
\label{eq:4.8} 
\overline{v}(x,t)\ge (e^{(t-\tau)\Delta}\overline{v}(\tau))(x)+(A-1)\int_\tau^t [e^{(t-s)\Delta}\overline{v}(s)^\frac{A}{A-1}](x)\, ds 
\end{equation} 
in $\mathbb{R}^N\times (\tau,T)$. 
Then, since 
$$
(A-1)
\int_\tau^t [e^{(t-s)\Delta}\overline{v}(s)^\frac{A}{A-1}](x)\, ds \le \overline{v}(x,t) 
$$
for all $\tau\in (0,t)$ and $\overline{v}$ is nonnegative, we see that
$$
\lim_{\tau\to 0}\int_\tau^t [e^{(t-s)\Delta}\overline{v}(s)^\frac{A}{A-1}](x)\, ds 
=\int_0^t [e^{(t-s)\Delta}\overline{v}(s)^\frac{A}{A-1}](x)\, ds 
$$
exists for all $(x,t)\in \mathbb{R}^N\times (0,T)$. 
Taking the limit $\tau\to 0$, by \eqref{eq:4.7} and \eqref{eq:4.8} we obtain 
$$
\overline{v}(x,t) 
\ge e^{t\Delta}v_0(x)+(A-1)\int_0^t [e^{(t-s)\Delta}\overline{v}(s)^\frac{A}{A-1}](x)\, ds, 
$$
that is, $\overline{v}$ is a supersolution of \eqref{eq:4.1}. 
Then we can construct a local in time 
solution $v\in C^{2,1}(\mathbb{R}^N\times (0,T))$ of the integral equation \eqref{eq:4.1} with the aid of Proposition~\ref{Proposition:2.1}. 
This yields a contradiction.

Next we consider the case $A=1$. 
The proof is similar to the above argument, thus we only show a brief sketch of the proof.
Let $r\in (0,N/2)$. 
By Proposition~\ref{Proposition:2.5} we can take a nonnegative function $w_0$ such that 
$e^{rw_0}\in {L}_{ul,\rho}^1(\mathbb{R}^N)$ and there can not exists a local in time 
solution for the 
integral 
equation 
\begin{equation} 
\label{eq:4.9} 
w(t,x)=e^{t\Delta}w_0+\int_0^t e^{(t-s)\Delta}e^{w(s)}ds. 
\end{equation} 
Define ${u}_0(x):=\max\{F^{-1}(e^{-w_0(x)}),s_2\}$. 
Suppose that there exists a local in time classical solution $\overline{u}\in C^{2,1}(\mathbb{R}^N\times (0,T))$ for problem 
\eqref{eq:1.1} with $\|\overline{u}(t)-e^{t\Delta}u_0\|_{L^{\infty}}\to 0$ as $t\to 0$.
As in the proof for the case $A>1$, 
$\overline{u}$
also satisfies 
the integral equation \eqref{eq:4.4}.
Define
$$
\overline{w}(x,t):=\log F(\overline u(x,t))^{-1} \quad\mbox{in}\,\,\, \mathbb{R}^N\times (0,T). 
$$
Similarly to the case $A>1$, we have $f'(\overline{u})F(\overline{u})\ge 1$ in $\mathbb{R}^N\times (0,T)$, 
and obtain  
$$
\partial_t \overline{w}-\Delta \overline{w}-e^{\overline{w}} 
=\frac{\left[f'(\overline u)F(\overline u)-1\right]|\nabla\overline u|^2}{f(\overline u)^2F(\overline{u})^2} \ge 0. 
$$
Thus, as in the above argument, by 
 the concavity of $\log F(s)^{-1}$ in $(s_2,\infty)$ 
 we see that $\overline{w}$ satisfies 
$$
\overline{w}(x,t)\ge (e^{t\Delta}w_0)(x)+\int_0^t [e^{(t-s)\Delta}e^{\overline{w}(s)}](x)\, ds, 
$$
and we can construct a solution of \eqref{eq:4.9} with the help of Proposition~\ref{Proposition:2.1}, 
which yields a contradiction. 
Thus we complete the proof of Theorem~\ref{Theorem:1.3}. 
\hfill 
\qed 

\vspace{5pt} 

\noindent 
{\bf Proof of Theorem~\ref{Theorem:1.4}.} 
Let $r\in (0, N/2)$. 
Put $h(s) := 1-f'(s)F(s)\ge 0$. 
Then we have $h(s) \to 0$ as $s\to\infty$ and 
$$
	(f(s)F(s))' = f'(s)F(s) -1 = -h(s) \le 0, 
$$
hence $l := \displaystyle\lim_{s\to\infty} f(s)F(s)\ge 0$ exists. 
Since $-h(s)$ is the derivative of $f(s)F(s)$ and $l$ exists, 
$h$ is integrable on $[s_3, \infty)$, 
where $s_3$ is the constant appearing in the assumption of Theorem~\ref{Theorem:1.4}. 
Let $\epsilon>0$ be a sufficiently small constant. 
For the case $l>0$, since 
$$
	(\log f(s))' = \frac{f'(s)F(s)}{f(s)F(s)} \ge \frac{1}{l}-\epsilon 
$$
for all sufficiently large $s>0$, 
by a simple calculation we get $f(s)\gtrsim e^{(1/l-\epsilon) s}$ for all sufficiently large $s>0$. 
On the other hand, for the case $l=0$, since 
\begin{equation} 
\label{eq:4.10} 
	f(s)F(s) = \int_s^\infty h(u)\, du =: \frac{1}{H(s)}, 
	\quad 
	(\log f(s))' 
	= 
	\frac{f'(s)F(s)}{f(s)F(s)} 
	\ge 
	(1-\epsilon)H(s), 
\end{equation} 
as in the above calculation, we obtain 
\begin{equation} 
\label{eq:4.11} 
	f(s) \gtrsim e^{(1-\epsilon) g(s)} 
	\quad\mbox{with} \,\,\, 
	g(s) = \int_{s_3}^s H(u)\, du 
\end{equation} 
for all $s\ge s_3$. 
Both cases $l=0$ and $l>0$ are treated in the same manner, 
so we give the proof only for the case $l=0$. 

Assume that, 
for any nonnegative initial function $u_0$ satisfying $F(u_0)^{-r}\in L_{ul, \rho}^1(\mathbb{R}^N)$, 
there exists a classical solution $u$ of 
$$
	u(t) = e^{t\Delta}u_0 + \int_0^t e^{(t-s)\Delta}f(u(s))\, ds. 
$$
Then, by \eqref{eq:4.11} we can take a constant $C>0$ such that 
\begin{equation} 
\label{eq:4.12} 
	u(t) \ge e^{t\Delta}u_0 + \int_0^t e^{(t-s)\Delta} e^{(1-\epsilon) g(u(s)) - C} \, ds. 
\end{equation} 
Remark that we can assume that $u_0\ge s_3$ without loss of generality 
by considering $\max\{u_0, s_3\}$ instead of $u_0$, and so $u\ge s_3$. 
Since 
$$
	g'(s) = H(s)>0, \qquad 
	\lim_{s\to\infty} \frac{g''(s)}{(g'(s))^2} 
	= \lim_{s\to\infty} h(s) 
	= 0, 
$$
we can apply the similar calculation as in the proof of Proposition~\ref{Proposition:2.5} for 
$$
	G(s) := \int_s^\infty \frac{du}{e^{(1-\epsilon)g(u) - C}}. 
$$
Then we obtain 
$$
	G(s) \gtrsim \frac{1}{g'(s)}e^{-(1-\epsilon) g(s)} 
	\iff 
	G(s)^{-r} \lesssim (g'(s))^r e^{(1-\epsilon)rg(s)} 
$$
for all sufficiently large $s>0$. 
On the other hand, by \eqref{eq:4.10} and \eqref{eq:4.11} we have 
$$
F(s)^{-r} \gtrsim (g'(s))^r e^{(1-\epsilon) r g(s)} 
$$
for all sufficiently large $s>0$. 
These imply that $G(s)^{-r}\lesssim F(s)^{-r}$ for all sufficiently large $s>0$. 
In particular, we have $G(u_0)^{-r} \in L^1_{ul, \rho}(\mathbb{R}^N)$. 
Then we see that \eqref{eq:4.12} yields a contradiction. 
In fact, if there exists a solution $u$ satisfying \eqref{eq:4.12}, 
then, in view of Proposition~\ref{Proposition:2.1}, we can construct a solution of 
$$
	u(t) = e^{t\Delta}u_0 + \int_0^t e^{(t-s)\Delta} e^{(1-\epsilon) g(u(s)) - C} \, ds. 
$$
This contradicts Proposition~\ref{Proposition:2.5} since $G(u_0)^{-r} \in L^1_{ul, \rho}(\mathbb{R}^N)$,  
and we complete the proof of Theorem~\ref{Theorem:1.4}. 
\hfill 
\qed

\section{Applications}
\label{section:5} 
In this section we apply Theorems~\ref{Theorem:1.1}--\ref{Theorem:1.4} 
to some examples of nonlinear heat equations. 
In particular, we consider the following cases: 
$f(u)=u^p+u^q$ ($p>q>1$) and $f(u)=e^{u^2}$. 
%
\subsection{Case $f(u)=u^p+u^q$ with $p>q>1$}  
Consider the case $f(u)=u^p+u^q$ with $p>q>1$, that is, 
\begin{equation} 
\label{eq:5.1} 
\left\{ 
\begin{array}{ll} 
\partial_t u=\Delta u+u^p+u^q, &x\in {\mathbb R}^N, \,\,\, t>0, \vspace{3pt} \\ 
u(x,0)=u_0(x)\ge 0, &x\in {\mathbb R}^N. 
\end{array} 
\right. 
\end{equation} 
%
Before stating the existence and nonexistence results for $f(u)=u^p+u^q$,
we prepare the following lemma. 
Recall that 
$$
F(s)=\int_{s}^\infty\frac{du}{u^p+u^q}. 
$$
\begin{Lemma} 
\label{Lemma:5.1}
Let $p>q>1$ and $f(s)=s^p+s^q$. 
Then there hold the following properties.
\begin{itemize}
\item[\rm (i)]
For all sufficiently large $s>0$, it holds 
\begin{equation}
\notag 
f'(s)F(s)\le \displaystyle \lim_{s\to \infty}f'(s)F(s)=\frac{p}{p-1}.
\end{equation}
\item[\rm (ii)] 
Let $r\ge 1$. 
Then 
\begin{equation} 
\notag 
\displaystyle F(s)^{-r} 
\lesssim s^{r(p-1)}+s^{r(q-1)}. 
\end{equation}
for all $s>0$. 
\end{itemize}
\end{Lemma}
%
%
\begin{Theorem}[$f(u)=u^p+u^q$]  
\label{Theorem:5.1} 
Let $N\ge 1$ and $p>q>1$. 
\begin{itemize} 
\item[\rm (i)] 
{\rm (Subcritical case)} 
Let $r>0$ satisfy $r\ge 1/(p-1)$ and $r>N/2$ and $u_0\in L^{r(p-1)}_{ul,\rho}(\mathbb{R}^N)$ be a nonnegative function.
Then there exists a local in time classical solution $u$ for problem~\eqref{eq:5.1} satisfying 
\begin{equation} 
\label{eq:5.2} 
\lim_{t\to 0}\|u(t)-e^{t\Delta}u_0\|_{L^{r(p-1)}_{ul,\rho}(\mathbb{R}^N)}=0. 
\end{equation} 
In particular, if $u_0\in \mathcal{L}^{r(p-1)}_{ul,\rho}(\mathbb{R}^N)$, then the solution converges to the initial data in $L^{r(p-1)}_{ul,\rho}(\mathbb{R}^N)$,
that is, 
$$
\lim_{t\to 0}\|u(t)-u_0\|_{L^{r(p-1)}_{ul,\rho}(\mathbb{R}^N)}=0.  
$$
Furthermore, the existence time $T$ can be estimated to satisfy 
\begin{equation} 
\label{eq:5.3} 
T^{\frac{N}{2p}}\rho^{-\frac{N}{p}}  
+ \max\left\{\|u_0\|_{L^{r(p-1)}_{ul,\rho}(\mathbb{R}^N)}^{r(p-1)}, \rho^N\right\} 
\left(T^{r-\frac{N}{2}\cdot\frac{p-1}{p}}\rho^{-N\cdot\frac{1}{p}}+T^{r-\frac{N}{2}}\right) 
\ge \gamma 
\end{equation}
for some $\gamma>0$ depending only on $N$, $p$, $q$ and $r$. 
\item[\rm (ii)] 
{\rm (Critical case)} 
Assume that $p>1+2/N$. 
Let $u_0\in \mathcal{L}^{\frac{N}{2}(p-1)}_{ul,\rho}(\mathbb{R}^N)$ be a nonnegative function. 
Then there exists a local in time classical solution $u$ for problem~\eqref{eq:5.1} satisfying 
\begin{equation} 
\label{eq:5.4} 
\lim_{t\to 0}\|u(t)-u_0\|_{L^{\frac{N}{2}(p-1)}_{ul,\rho}(\mathbb{R}^N)}=0. 
\end{equation} 
\item[\rm (iii)] 
{\rm (Nonexistence)} 
Let $p>1+2/N$ and $1/(p-1)\le r<N/2$. 
Then there exists a nonnegative initial function $u_0\in {L}^{r(p-1)}_{ul,\rho}(\mathbb{R}^N)$ 
and problem~\eqref{eq:5.1} can not possess local in time nonnegative classical  
solutions satisfying \eqref{eq:5.2}. 
\end{itemize} 
\end{Theorem} 
For the proof of Theorem~\ref{Theorem:5.1},
we check the conditions of Theorems~\ref{Theorem:1.1} and \ref{Theorem:1.3}. 
To this end, we start from the proof of Lemma~\ref{Lemma:5.1}.
\vspace{5pt} 

\noindent 
{\bf Proof of Lemma~\ref{Lemma:5.1}.} 
We see that $f$ is a positive convex function in $(0,\infty)$ and 
\begin{align*} 
f'(s)F(s) 
&= (ps^{p-1}+qs^{q-1}) \int_s^\infty \frac{du}{u^p+u^q} \\ 
&= (ps^{p-1}+qs^{q-1}) \int_s^\infty \frac{du}{u^p(1+u^{q-p})} \\ 
&= (ps^{p-1}+qs^{q-1}) \left[ 
\frac{s^{1-p}}{(p-1)(1+s^{q-p})}+\frac{p-q}{p-1}\int_s^\infty u^{-(p-1)}\cdot \frac{u^{q-p-1}}{(1+u^{q-p})^2}\, du 
\right] \\ 
&\le \frac{p+qs^{q-p}+o(s^{q-p})}{(p-1)(1+s^{q-p})} 
<\frac{p+ps^{q-p}}{(p-1)(1+s^{q-p})} 
=\frac{p}{p-1} 
\end{align*} 
for all sufficiently large $s>0$. 
Here we used 
$$
\int_s^\infty u^{-(p-1)}\cdot \frac{u^{q-p-1}}{(1+u^{q-p})^2}\, du\le \int_s^\infty u^{q-p-1-(p-1)}\, du 
=O(s^{q-p-(p-1)})  
=o(s^{q-p})
$$
as $s\to\infty$. 
Furthermore, it is easy to check that 
$$
\lim_{s\to\infty}f'(s)F(s)=\frac{p}{p-1}=1+\frac{1}{p-1}>1. 
$$
This proves (i).


It remains to prove (ii).
Since 
\begin{equation} 
\notag 
\begin{split}
F(s) &= \int_s^\infty \frac{du}{u^q(1+u^{p-q})} 
=\frac{s^{1-q}}{(q-1)(1+s^{p-q})}-\frac{p-q}{q-1}\int_s^\infty u^{1-q}\cdot \frac{u^{p-q-1}}{(1+u^{p-q})^2}\, du \\ 
&= \frac{s^{1-q}}{q-1}\cdot \frac{1}{1+s^{p-q}}-\frac{p-q}{q-1}\int_s^\infty \frac{1}{u^p+u^q}\cdot \frac{u^{p-q}}{1+u^{p-q}}\, du, 
\end{split}
\end{equation}  
we have 
\begin{equation} 
\label{eq:5.8} 
F(s)\ge \frac{s^{1-q}}{2(q-1)}-\frac{p-q}{q-1}F(s) \iff F(s)\ge \frac{1}{2(p-1)}s^{1-q} 
\end{equation} 
for all $s\in (0,1)$. 
On the other hand, 
since $f(s)=s^p+s^q\le 2s^p$ for $s\ge 1$, we have 
\begin{equation} 
\label{eq:5.9} 
F(s)\ge \frac{1}{2}\int_s^\infty \frac{du}{u^p} =\frac{s^{-(p-1)}}{2(p-1)} 
\end{equation} 
for all $s\ge 1$. 
Combining \eqref{eq:5.8} and \eqref{eq:5.9}, we obtain $F(s)^{-r}\lesssim s^{r(p-1)}+s^{r(q-1)}$ for all $s>0$.  
This yields the assertion~(ii). 
Thus we complete the proof of Lemma~\ref{Lemma:5.1}. 
\hfill \qed 

\vspace{5pt} 

\noindent 
{\bf Proof of Theorem~\ref{Theorem:5.1}.} 
Applying Lemma~\ref{Lemma:5.1},
we check the conditions of Theorems~\ref{Theorem:1.1}--\ref{Theorem:1.3}. 
We first prove assertion~(i). 
By Lemma~\ref{Lemma:5.1}~(ii) we have 
\begin{equation}
\label{eq:a}  
\begin{split} 
\int_{B_{\rho}(y)}F(u_0(x))^{-r}\, dx 
&= \left(\int_{B_{\rho}(y)\cap \{u_0(x)<1\}}+\int_{B_{\rho}(y)\cap \{u_0(x)\ge 1\}}\right)F(u_0(x))^{-r}\, dx \\ 
&\le F(1)^{-r}|B_{\rho}(0)|+(2(p-1))^{r}\int_{B_{\rho}(y)}|u_0(x)|^{r(p-1)}\, dx 
\end{split} 
\end{equation} 
for all $y\in {\mathbb R}^N$. 
Then, since $u_0\in L^{r(p-1)}_{ul,\rho}(\mathbb{R}^N)$, 
we have $F(u_0)^{-r}\in {L}^1_{ul,\rho}(\mathbb{R}^N)$. 
By Lemma~\ref{Lemma:5.1}~(i) we have 
$$ 
r\ge \frac{1}{p-1}=\lim_{s\to\infty}f'(s)F(s)-1, 
$$ 
and we see from Theorem~\ref{Theorem:1.1}~(i) that 
problem~\eqref{eq:5.1} has a local in time solution $u$ satisfying \eqref{eq:5.2}. 
Furthermore, by \eqref{eq:1.13} and \eqref{eq:a} we obtain the estimate on the existence time and prove \eqref{eq:5.3}.  

We next prove assertion~(ii). 
Since $p>1+2/N$, by Lemma~\ref{Lemma:5.1}~(i) we have 
$$
\frac{N}{2}>\frac{1}{p-1}=\lim_{s\to\infty}f'(s)F(s)-1. 
$$
Assuming that $u_0\in {\mathcal{L}^{\frac{N}{2}(p-1)}_{ul,\rho}(\mathbb{R}^N)}$, 
we prove that the initial data $u_0$ satisfies
\begin{equation} 
\label{eq:5.10} 
F(u_0)^{-\frac{N}{2}}\in {\cal L}_{ul,\rho}^{1}(\mathbb{R}^N).
\end{equation} 
Since 
$u_0\in {\mathcal{L}^{\frac{N}{2}(p-1)}_{ul,\rho}(\mathbb{R}^N)}$,
there exists a sequence
$\{u_n\}\subset BUC(\mathbb{R}^N)$
such that
$u_n\to u_0$ in $ {{L}^{\frac{N}{2}(p-1)}_{ul,\rho}(\mathbb{R}^N)}$ as $n\to \infty$. 
Then, by the H\"older inequality we have 
$u_n\to u_0$ in $ {{L}^{\frac{N}{2}(q-1)}_{ul,\rho}(\mathbb{R}^N)}$.
Remark that 
$F(u_n)^{-\frac{N}{2}}\in BUC(\mathbb{R}^N)$.
Applying the mean value theorem, we have
\[
 \left|
  F(u_0)^{-\frac{N}{2}}- F(u_n)^{-\frac{N}{2}}
 \right|
 \le  \frac{N}{2} \frac{\displaystyle F(u_0+\theta(u_n-u_0))^{-\frac{N}{2}-1}}{\displaystyle f(u_0+\theta(u_n-u_0))}
|u_0-u_n|, 
\]
where $0<\theta<1$.
Let 
$v:=u_0+\theta(u_n-u_0)$.
Since $f(s)\ge s^p$ and $f(s)\ge s^q$, by Lemma~\ref{Lemma:5.1}~{(ii)} we have  
\[
 \left|
  F(u_0)^{-\frac{N}{2}}- F(u_n)^{-\frac{N}{2}}
 \right|
\lesssim \left(v^{\frac{N}{2}(p-1)-1}+v^{\frac{N}{2}(q-1)-1}\right) |u_0-u_n|.
\]
Therefore it follows from the H\"older inequality that 
\[
\begin{split}
\int_{B_{\rho}(y)}
 \left|
  F(u_0)^{-\frac{N}{2}}- F(u_n)^{-\frac{N}{2}}
 \right|
 dx 
&
\lesssim 
\int_{B_{\rho}(y)}
\left(v^{\frac{N}{2}(p-1)-1}+v^{\frac{N}{2}(q-1)-1}\right) |u_0-u_n| \, dx
\\
&
\lesssim 
\left(
\int_{B_{\rho}(y)}
v^{\frac{N}{2}(p-1)}
dx
\right)^{\frac{p-1-\frac2N}{(p-1)}} 
\|u_0-u_n\|_{L^{\frac{N}{2}(p-1)}(B_{\rho}(y))}
\\
&\,\,\,\, 
+
\left(
\int_{B_{\rho}(y)}
v^{\frac{N}{2}(q-1)}
dx
\right)^{\frac{q-1-\frac2N}{(q-1)}} 
\|u_0-u_n\|_{L^{\frac{N}{2}(q-1)}(B_{\rho}(y))}.
\end{split}
\]
Then $F(u_n)^{-\frac{N}{2}}$  
converges to $F(u_0)^{-\frac{N}{2}}$ as $n\to\infty$ in $L^1_{ul,\rho}(\mathbb{R}^N)$ 
since
$$
u_0,u_n\in L^{\frac{N}{2}(p-1)}_{ul,\rho}(\mathbb{R}^N)\subset L^{\frac{N}{2}(q-1)}_{ul,\rho}(\mathbb{R}^N)\ \ 
\text{and}\ u_n\to u_0\ \ \text{in}\ 
L^{\frac{N}{2}(p-1)}_{ul,\rho}(\mathbb{R}^N),  
$$
and so we obtain \eqref{eq:5.10}. 
Therefore, by Theorem~\ref{Theorem:1.1}~(ii) 
we get a local in time solution for problem~\eqref{eq:5.1} satisfying \eqref{eq:5.2} with $r=N/2$. 
Since $u_0\in \mathcal{L}^{\frac{N}{2}(p-1)}_{ul,\rho}(\mathbb{R}^N)$, 
by Lemma~\ref{Lemma:2.2} we see that $e^{t\Delta}u_0$ converges to $u_0$ in ${L^{\frac{N}{2}(p-1)}_{ul,\rho}(\mathbb{R}^N)}$ as $t\to 0$. Thus
\eqref{eq:5.4} holds.


We finally prove assertion~(iii). 
For the proof of assertion~(iii), we apply Theorem~\ref{Theorem:1.3}. 
Let $p>1+2/N$. 
Consider 
\begin{equation} 
\label{eq:5.11} 
\left\{ 
\begin{array}{ll} 
\partial_t u=\Delta u+u^p, &x\in \mathbb{R}^N, \,\,\, t>0, \vspace{3pt} \\ 
u(x,0)=u_0(x)\ge 0, &x\in \mathbb{R}^N. 
\end{array}  
\right. 
\end{equation} 
Put $f_0(s)=s^p$. 
Then, since 
$$
f_0'(s)\int_s^\infty \frac{du}{f_0(u)}=\frac{p}{p-1}, \qquad \frac{p}{p-1}-1=\frac{1}{p-1}<\frac{N}{2}, 
$$
for $r\in [1/(p-1), N/2)$, by Theorem~\ref{Theorem:1.3} we find an initial function 
$u_0\in {L}_{ul,\rho}^{r(p-1)}(\mathbb{R}^N)$ 
such that 
there can not exist nonnegative classical solutions of \eqref{eq:5.11}. 
See also \cite[Corollary~5.1]{W1}. 
Suppose that there exists a classical solution $u$ of \eqref{eq:5.1} with this initial data $u_0$. 
Then we have 
$$
\partial_t u\ge \Delta u+u^p, \quad x\in\mathbb{R}^N, \,\,\, t>0, \qquad u(x,0)=u_0(x), \quad x\in\mathbb{R}^N, 
$$
and so $u$ is a supersolution of \eqref{eq:5.11}. 
Thus, as in the proof of Theorem~\ref{Theorem:1.3}, by Proposition~\ref{Proposition:2.1} we can construct a solution of \eqref{eq:5.11}, 
which is a contradiction. 
Therefore we prove assertion~(iii), and the proof of Theorem~\ref{Theorem:5.1} is complete. 
\hfill 
\qed 

\subsection{Case $f(u)=e^{u^2}$} 
Consider the case $f(u)=e^{u^2}$, that is, 
\begin{equation} 
\label{eq:5.12} 
\left\{ 
\begin{array}{ll} 
\partial_t u=\Delta u+e^{u^2}, &x\in {\mathbb R}^N, \,\,\, t>0, \vspace{3pt} \\ 
u(x,0)=u_0(x)\ge 0, &x\in {\mathbb R}^N. 
\end{array} 
\right. 
\end{equation} 
Recall that 
\begin{equation} 
\label{eq:5.13} 
F(s)=\int_s^\infty\frac{du}{e^{u^2}}. 
\end{equation} 
We first 
prepare several lemmas. 
%
\begin{Lemma} 
\label{Lemma:5.2} 
Let $f(s)=e^{s^2}$ and $F$ be the function defined by \eqref{eq:5.13}. 
\begin{itemize} 
\item[\rm (i)] 
It holds $F(s)^{-1}\lesssim (1+s)e^{s^2}$ for all $s>0$. 
\item[\rm (ii)] 
Let $0<\sigma \le 1$.
Define
$h_{\sigma}(t):=t^{\sigma} e^{\sigma t^2}$ and $g_{\sigma}(s):=F(h_{\sigma}^{-1}(s))^{-\sigma}$, 
where $h_\sigma^{-1}$ denotes the inverse function of $h_\sigma$. 
Then there holds that  
$|g_{\sigma}'(s)|\lesssim 1$ for all $s>0$. 
\end{itemize} 
\end{Lemma} 
\noindent 
{\bf Proof.} 
By integration by parts we have 
\begin{equation} 
\notag 
F(s)=\frac{1}{2}\int_s^\infty u^{-1}(-e^{-u^2})'\, du=\frac{1}{2}s^{-1}e^{-s^2}-\frac{1}{2}\int_s^\infty u^{-2}e^{-u^2}\, du 
\end{equation} 
for all $s>0$. 
Since 
\begin{equation} 
\notag 
\int_s^\infty u^{-2}e^{-u^2}\, du = \int_s^\infty u^{-3}\cdot ue^{-u^2}\, du 
\le s^{-3}\int_s^\infty ue^{-u^2}\, du 
=\frac{1}{2}s^{-3}e^{-s^2} 
\le \frac{1}{2}s^{-1}e^{-s^2} 
\end{equation} 
for all $s\ge 1$, 
we obtain 
\begin{equation} 
\label{eq:5.14} 
F(s)\ge \frac{1}{4}s^{-1}e^{-s^2} 
\end{equation} 
for all $s\ge 1$. 
In particular, since $F(0)>0$, we have $F(0)^{-1}<\infty$, 
and by \eqref{eq:5.14} we obtain assertion~(i). 

We next prove assertion~(ii). 
Put $t(s):=h_{\sigma}^{-1}(s)$. 
By the definition of $h_{\sigma}$, we have 
$s=t(s)^{\sigma} e^{\sigma t(s)^2}$.  
Thus we obtain 
$$ 
1=\left(\sigma +2 \sigma t(s)^{2}\right)t(s)^{\sigma -1}e^{\sigma t(s)^2}t'(s).
$$
This together with assertion~(i) and the assumption $0<\sigma\le 1$ implies that 
\begin{align*} 
|g_{\sigma}'(s)| &= \sigma F(t)^{-\sigma-1}f(t)^{-1} t'(s) 
\lesssim (1+t(s))^{\sigma+1}e^{\sigma t(s)^2}\cdot \frac{t(s)^{1-\sigma}e^{-\sigma t(s)^2}}
{\sigma+2\sigma t(s)^2} 
\lesssim 1  
\end{align*} 
for all $s>0$.  
This proves assertion~(ii). 
Thus we complete the proof of Lemma~\ref{Lemma:5.2}. 
\hfill \qed 

\vspace{5pt} 

With the help of Lemma~\ref{Lemma:5.2}, we can prove the following assertion. 
\begin{Lemma} 
\label{Lemma:5.3} 
Let $r>0$ and $h_r(u_0)=|u_0|^r e^{r|u_0(x)|^2}\in {\cal L}_{ul,\rho}^1(\mathbb{R}^N)$. 
Then $F(u_0)^{-r}\in {\cal L}_{ul,\rho}^1(\mathbb{R}^N)$. 
\end{Lemma} 
\noindent 
{\bf Proof.} 
Let $r>0$ and assume $h_r(u_0)=|u_0|^r e^{r|u_0(x)|^2}\in {\cal L}_{ul,\rho}^1(\mathbb{R}^N)$.
Put $R:= \max\{r,1\}$ and $\sigma:=\frac{r}{R}$, then clearly $0<\sigma \le 1$.
Recall that $h_r(u_0)=|u_0|^r e^{r|u_0(x)|^2}\in {\cal L}_{ul,\rho}^1(\mathbb{R}^N)$ 
is equivalent to $h_{\sigma}(u_0)=|u_0|^{\sigma} e^{\sigma |u_0(x)|^2}\in {\cal L}_{ul,\rho}^R(\mathbb{R}^N)$.
Now define $v_0(x):=h_{\sigma}(u_0)$, 
then $v_0\in {\cal L}_{ul,\rho}^R(\mathbb{R}^N)$
and $g_{\sigma}(v_0(x))=F(u_0(x))^{-\sigma}$, 
where $g_{\sigma}$ is the function defined in Lemma~\ref{Lemma:5.2}. 
Since $v_0\in {\cal L}_{ul,\rho}^R(\mathbb{R}^N)$,  one can take  a sequence $\{v_n\}_{n\in\mathbb{N}}\subset BUC(\mathbb{R}^N)$ such that 
$v_n\to v_0$ in $L_{ul,\rho}^R(\mathbb{R}^N)$ as $n\to\infty$. 
Define $u_n:=g_{\sigma}(v_n)$, then $u_n$ also belongs to $BUC(\mathbb{R}^N)$.
This together with the mean value theorem and Lemma~\ref{Lemma:5.2}~(ii) implies that 
\begin{align*} 
\sup_{y\in \mathbb{R}^N}\int_{B_{\rho}(y)}|u_n(x)-F(u_0(x))^{-\sigma}|^{R}\, dx 
&=\sup_{y\in \mathbb{R}^N}\int_{B_{\rho}(y)}|g_\sigma(v_n(x))-g_\sigma(v_0(x))|^R\, dx \\ 
&\lesssim 
\sup_{y\in\mathbb{R}^N}\int_{B_{\rho}(y)}|v_n(x)-v_0(x)|^R\, dx \to 0 
\end{align*} 
as $n\to\infty$. 
Thus we see that $F(u_0)^{-\sigma}\in {\cal L}^R_{ul,\rho}(\mathbb{R}^N)$, 
which implies
$F(u_0)^{-r}\in {\cal L}^1_{ul,\rho}(\mathbb{R}^N)$. 
Therefore we complete the proof of Lemma~\ref{Lemma:5.3}. 
\hfill \qed 

\vspace{5pt} 

Furthermore, we introduce one lemma on the convergence of $e^{t\Delta}u_0$ to $u_0$ as $t\to 0$. 
\begin{Lemma} 
\label{Lemma:5.4} 
Let $r>0$ and $h_r(u_0)=|u_0|^r e^{r|u_0|^2}\in \mathcal{L}_{ul,\rho}^1(\mathbb{R}^N)$, then 
$$ 
\lim_{t\to 0}\|h_r(|e^{t\Delta}u_0-u_0|)\|_{L^1_{ul,\rho}(\mathbb{R}^N)}=0. 
$$ 
\end{Lemma} 
\noindent 
{\bf Proof.} 
By the assumption we have $h_r(u_0)\in {\cal L}_{ul,\rho}^1(\mathbb{R}^N)$. 
Note that $u_0$ is a nonnegative function. 
By Lemma~\ref{Lemma:2.2} we have 
\begin{equation} 
\label{eq:5.15} 
\sup_{y\in\mathbb{R}^N}\int_{B_{\rho}(y)}|h_r(u_0(x+z))-h_r(u_0(x))|\, dx\to 0 
\end{equation} 
as $z\to 0$. 
Then we prove that 
\begin{equation} 
\label{eq:5.16} 
\sup_{y\in \mathbb{R}^N}\int_{B_{\rho}(y)}h_r(|u_0(x+z)-u_0(x)|)\, dx\to 0 
\end{equation} 
as $z\to 0$. 
For the proof of \eqref{eq:5.16}, we first assume that $r\ge 1$. 
Let $s\ge t\ge 0$. 
Then, 
by an elementary inequality 
\begin{equation} \notag
|s-t|^p\le |s^p-t^p| \quad\mbox{for}\,\,\, s,t\ge 0, \ \ p\ge 1,
\end{equation} 
%
we have 
$$
h_r(|s-t|)=\sum_{n=0}^\infty \frac{r^n}{n!} (s-t)^{2n+r} 
\le \sum_{n=0}^\infty \frac{r^n}{n!} (s^{2n+r}-t^{2n+r})=h_r(s)-h_r(t). 
$$
Similarly, we obtain $h_r(|s-t|)\le h_r(t)-h_r(s)$ for $t\ge s\ge 0$. 
Thus we have $h_r(|s-t|)\le |h_r(s)-h_r(t)|$ for all $s$, $t\ge 0$. 
This together with \eqref{eq:5.15} gives \eqref{eq:5.16}. 
Next we consider the case $0<r<1$. 
Put $\tilde{h}(s):=h_r(s)-s^r=s^r(e^{rs^2}-1)$ and $\hat{h}(s):=h_r(s^{1/r})=se^{rs^{2/r}}$. 
Then we have $\hat{h}'(s)=e^{rs^{2/r}}+2s^{2/r}e^{rs^{2/r}}\ge 1$, and by the mean value theorem we obtain 
$|\hat{h}(s)-\hat{h}(t)|=|\hat{h}'(\theta)||s-t|\ge |s-t|$ for all $s$, $t\ge 0$, where $\theta\in (0,1)$. 
This is equivalent to $|s^r-t^r|\le |h_r(s)-h_r(t)|$. 
Then we have 
\begin{equation} 
\label{eq:5.17} 
\int_{B_{\rho}(y)}|u_0(x+z)^r-u_0(x)^r|\, dx \le \int_{B_{\rho}(y)}|h_r(u_0(x+z))-h_r(u_0(x))|\, dx 
\end{equation} 
for all $y$, $z\in \mathbb{R}^N$. 
On the other hand, similarly to the calculation for $r\ge 1$, 
we see that $\tilde{h}(|s-t|)\le |\tilde{h}(s)-\tilde{h}(t)|$ for all $s$, $t\ge 0$. 
Therefore, since $\tilde{h}(s)=h_r(s)-s^r$, by \eqref{eq:5.17} we obtain 
\begin{equation} 
\label{eq:5.18} 
\begin{split} 
&\int_{B_{\rho}(y)}\tilde{h}(|u_0(x+z)-u_0(x)|)\, dx 
\le \int_{B_{\rho}(y)}|\tilde{h}(u_0(x+z))-\tilde{h}(u_0(x))|\, dx \\ 
&\qquad \le \int_{B_{\rho}(y)}|h_r(u_0(x+z))-h_r(u_0(x))|\, dx+\int_{B_{\rho}(y)}|u_0(x+z)^r-u_0(x)^r|\, dx \\ 
&\qquad \le 2\int_{B_{\rho}(y)}|h_r(u_0(x+z))-h_r(u_0(x))|\, dx. 
\end{split} 
\end{equation}  
Furthermore, 
by the H\"{o}lder inequality we obtain 
$$
\int_{B_{\rho}(y)}|u(x)-v(x)|^r\, dx 
\le |B_{\rho}(0)|^\frac{2}{r+2}\left(\int_{B_{\rho}(y)}|u(x)-v(x)|^{r+2}\, dx\right)^\frac{r}{r+2} 
$$
for all $y\in \mathbb{R}^N$ and suitable measurable functions $u$, $v$. 
This together with $\tilde{h}(s)\ge rs^{r+2}$ implies that 
\begin{equation} 
\label{eq:5.19} 
\int_{B_{\rho}(y)}|u(x)-v(x)|^r\, dx \le |B_{\rho}(0)|^\frac{2}{r+2}\left(\frac{1}{r}\int_{B_{\rho}(y)}\tilde{h}(|u(x)-v(x)|)\, dx\right)^\frac{r}{r+2}. 
\end{equation} 
Then, since $h_r(s)=\tilde{h}(s)+s^r$, 
by \eqref{eq:5.15}, \eqref{eq:5.18} and \eqref{eq:5.19} with $u=u_0(\cdot+z)$, $v=u_0$ we obtain \eqref{eq:5.16}. 

Once we get \eqref{eq:5.16}, we can easily prove the lemma. 
We give the proof only for the case $0<r<1$. 
Since $\tilde{h}$ is a convex function and 
$$
(e^{t\Delta}u_0)(x)-u_0(x)=(4\pi)^{-\frac{N}{2}}\int_{\mathbb{R}^N}e^{-\frac{|w|^2}{4}}(u_0(x+\sqrt{t}w)-u_0(x))\, dw, 
$$
by the Jensen inequality and the Fubini theorem we have 
\begin{equation} 
\label{eq:5.20} 
\int_{B_{\rho}(y)}\tilde{h}(|e^{t\Delta}u_0-u_0|)\, dx 
\le (4\pi)^{-\frac{N}{2}}\int_{\mathbb{R}^N}e^{-\frac{|w|^2}{4}}\int_{B_{\rho}(y)}\tilde{h}(|u_0(x+\sqrt{t}w)-u_0(x)|)\, dx\, dw 
\end{equation} 
for all $y\in \mathbb{R}^N$ and $t>0$. 
Furthermore, by \eqref{eq:5.19} with $u=e^{t\Delta}u_0$ and $v=u_0$ we have 
\begin{equation} 
\label{eq:5.21} 
\int_{B_{\rho}(y)}|e^{t\Delta}u_0-u_0|^r\, dx\le |B_{\rho}(0)|^\frac{2}{r+2}\left( 
\frac{1}{r}\int_{B_{\rho}(y)}\tilde{h}(|e^{t\Delta}u_0-u_0|)\, dx  
\right)^\frac{r}{r+2}. 
\end{equation} 
Since $h_r(s)=\tilde{h}(s)+s^r$ and $\tilde{h}(s)\le h_r(s)$, 
by \eqref{eq:5.16}, \eqref{eq:5.20} and \eqref{eq:5.21} we obtain the desired convergence and complete the proof of Lemma~\ref{Lemma:5.4}. 
\hfill \qed 

\vspace{5pt} 

We are ready to state the results on existence and nonexistence of solutions for problem~\eqref{eq:5.12}. 
\begin{Theorem}[$f(u)=e^{u^2}$] 
\label{Theorem:5.2} 
Let $N\ge 1$. 
\begin{itemize}  
\item[\rm (i)] 
{\rm (Subcritical case)} 
For any $r>N/2$ and nonnegative initial function $u_0$ satisfying 
$$ 
  |u_0|^r e^{r|u_0|^2}\in L_{ul,\rho}^1(\mathbb{R}^N), 
$$ 
there exists a local in time classical solution $u$ for problem~\eqref{eq:5.12} such that 
\begin{equation} 
\label{eq:5.22} 
\lim_{t\to 0}\|u(t)-e^{t\Delta}u_0\|_{L^\infty(\mathbb{R}^N)}=0. 
\end{equation} 
In addition, if $|u_0|^r e^{r|u_0|^2}\in {\cal L}_{ul,\rho}^1(\mathbb{R}^N)$, then 
\begin{equation} 
\label{eq:5.23} 
\lim_{t\to 0}\, \sup_{y\in\mathbb{R}^N}\int_{B_{\rho}(y)}|u(x,t)-u_0(x)|^r e^{r|u(x,t)-u_0(x)|^2}\, dx=0. 
\end{equation} 
Furthermore, for any sufficiently small $\epsilon>0$, the existence time $T$ can be chosen to satisfy
\begin{equation} 
\label{eq:5.24}
T^{\frac{\epsilon N}{2}} \rho^{-\epsilon N} 
+\gamma_\epsilon \max\left\{\|u_0^r e^{ru_0^2}\|_{L_{ul,\rho}^1(\mathbb{R}^N)}, \rho^N\right\} 
\left(T^{r-\frac{N}{2}(1-\epsilon)}\rho^{-\epsilon N}+T^{r-\frac{N}{2}}\right) 
\ge 
\gamma 
\end{equation}
where $\gamma$ depends only on $N$ and $r$, 
and $\gamma_\epsilon$ is a positive constant depending only on $N$, $r$ and $\epsilon$ 
satisfying $\gamma_\epsilon\to \infty $ as $\epsilon\to 0$. 
\item[\rm (ii)] 
{\rm (Critical case)} 
Let $u_0$ be a nonnegative initial function such that $|u_0|^\frac{N}{2} e^{\frac{N}{2}|u_0|^2}\in {\cal L}_{ul,\rho}^1(\mathbb{R}^N)$. 
Then there exists a local in time classical solution for problem~\eqref{eq:5.12} 
satisfying \eqref{eq:5.22} and \eqref{eq:5.23} with $r=N/2$. 
\item[\rm (iii)] 
{\rm (Nonexistence)} 
Let $0<r<N/2$. 
Then there exists a nonnegative initial function $u_0$ such that $|u_0|^r e^{r|u_0|^2}\in L_{ul,\rho}^1(\mathbb{R}^N)$ 
and problem~\eqref{eq:5.12} with \eqref{eq:5.22}  
can not possess local in time nonnegative classical solutions. 
\end{itemize} 
\end{Theorem} 
\noindent 
{\bf Proof.} 
We first prove assertion~(i). 
Let $f(u)=e^{u^2}$. 
Then we have  
\begin{equation} 
\label{eq:5.25} 
f'(s)F(s)=2se^{s^2}\int_s^\infty \frac{du}{e^{u^2}}\le 
e^{s^2}\int_s^\infty 2ue^{-u^2}\, du=1 
\end{equation} 
for all $s>0$. 
Furthermore, we have $\displaystyle\lim_{s\to\infty}f'(s)F(s)=1$. 
For $r>N/2$, let $u_0$ satisfy $|u_0|^r e^{r|u_0|^2}\in L_{ul,\rho}^1(\mathbb{R}^N)$. 
Then, by \eqref{eq:5.14} we have 
\begin{align*} 
\sup_{y\in {\mathbb R}^N} \int_{B_{\rho}(y)}F(u_0(x))^{-r}\, dx 
&\le F(1)^{-r}|B_{\rho}(0)| +4^r \sup_{y\in {\mathbb R}^N} \int_{B_{\rho}(y)}|u_0(x)|^r e^{r|u_0(x)|^2}\, dx<\infty. 
\end{align*} 
Then we can apply Theorem~\ref{Theorem:1.1}~(i) to problem~\eqref{eq:5.12}, 
and obtain local in time existence of a classical solution for \eqref{eq:5.12} satisfying \eqref{eq:5.22}. 
The convergence \eqref{eq:5.23} follows directly from Lemma~\ref{Lemma:5.4}. 
The estimate of the existence time \eqref{eq:5.24} follows from \eqref{eq:1.15}. 

We next prove assertion~(ii).  
If $|u_0|^\frac{N}{2}e^{\frac{N}{2}|u_0|^2}\in {\cal L}_{ul,\rho}^1(\mathbb{R}^N)$, 
then by Lemma~\ref{Lemma:5.3} we have 
$F(u_0)^{-\frac{N}{2}}\in {\cal L}_{ul,\rho}^1(\mathbb{R}^N)$. 
Thus, by Theorem~\ref{Theorem:1.1}~(ii) we can obtain a local in time classical solution for \eqref{eq:5.12} satisfying \eqref{eq:5.22}. 
We also obtain \eqref{eq:5.23} with $r=N/2$ by Lemma~\ref{Lemma:5.4}. 

We finally prove assertion~(iii). 
Let $r\in (0, N/2)$. 
By \eqref{eq:5.25} we have $f'(s)F(s) \le 1$ for all $s>0$. 
Then, by Theorem~\ref{Theorem:1.4} we can take a initial function $u_0$ satisfying 
$F(u_0)^{-r}\in L^1_{ul, \rho}(\mathbb{R}^N)$ such that 
there can not exist a solution for problem~\eqref{eq:5.12} with \eqref{eq:5.22}. 
On the other hand, by \eqref{eq:5.25} we have  
$F(s)^{-r}\ge 2^r s^r e^{rs^2}$
for all $s>0$. 
This implies that $|u_0|^r e^{r|u_0|^2}\in L_{ul, \rho}^1(\mathbb{R}^N)$, 
and we complete the proof of Theorem~\ref{Theorem:5.2}. 
\hfill 
\qed 

\appendix 
\section{Local in time existence of solutions for a semilinear heat equation in the uniformly local $L^r$ spaces} 
\label{section:A}
Consider a semilinear heat equation 
\begin{equation} 
\label{eq:A.1}%
\left\{ 
\begin{array}{ll} 
\partial_t u=\Delta u+|u|^{p-1}u, &x\in {\mathbb{R}}^N, \,\,\, t>0, \vspace{3pt} \\ 
u(x,0)=u_0(x), &x\in {\mathbb{R}}^N, 
\end{array} 
\right. 
\end{equation} 
where $p>1$. 
In this section we consider the case $u_0$ belongs to a uniformly local $L^r$ space, and 
study existence of local in time solutions of the integral equation 
$$
u(t)=e^{t\Delta}u_0+\int_0^t e^{(t-s)\Delta}|u(s)|^{p-1}u(s)\, ds. 
$$
In particular, we prove Propositions~\ref{Proposition:2.2} and \ref{Proposition:2.3}. 
These propositions can be proved by the similar argument as in \cite{W1} 
with the aid of Lemma~\ref{Lemma:2.1}. 
More precisely, we prove that the map $\Phi$ defined by 
$$
  \Phi(u)
  :=e^{t\Delta}u_0+\int_0^t e^{(t-s)\Delta}|u(s)|^{p-1}u(s)\, ds, 
$$
is a contraction map from a suitable Banach space to itself. 
\vspace{5pt} 

\noindent 
{\bf Proof of Proposition~\ref{Proposition:2.2}.} 
Let $r\ge 1$, $r>\frac{N}{2}(p-1)$ and $u_0\in L^r_{ul,\rho}(\mathbb{R}^N)$.
Define 
$$
  X_{M,T}
  :=
  \left\{  
  \begin{array}{l} 
   u\in L^{\infty}(0,T; L^{r}_{ul,\rho}(\mathbb{R}^N))
   \cap L^{\infty}_{loc}((0,T); L^{pr}_{ul,\rho}(\mathbb{R}^N))
   :\, \vspace{5pt} \\ 
   \displaystyle\sup_{0<t<T}\|u\|_{L^r_{ul,\rho}(\mathbb{R}^N)}\le M, 
   \displaystyle\sup_{0<t<T}t^{\alpha}\|u(t)\|_{L^{pr}_{ul,\rho}(\mathbb{R}^N)}\le M
  \end{array}
  \right\},
$$ 
equipped with the metric $d_X(u,v):=\sup\{t^\alpha \|u(t)-v(t)\|_{L_{ul,\rho}^{pr}(\mathbb{R}^N)}:\, 0<t<T\}$, 
where $\alpha:=\frac{N}{2}\left(\frac{1}{r}-\frac{1}{pr}\right)$ 
and $M$, $T$ are positive constants to be chosen later. 
Remark that $\alpha p<1$. 
For $u\in X_{M,T}$, by Lemma~\ref{Lemma:2.1} we have 
\begin{equation} 
\notag 
 \begin{split}
  \|\Phi(u)\|_{L^r_{ul,\rho}(\mathbb{R}^N)}
  &
  \le 
  \|u_0\|_{L^r_{ul,\rho}(\mathbb{R}^N)}+\int_0^t \|u(s)\|_{L^{pr}_{ul,\rho}(\mathbb{R}^N)}^p\, ds
  \\
  &
  \le 
  \|u_0\|_{L^r_{ul,\rho}(\mathbb{R}^N)}+\int_0^t s^{-\alpha p}\, ds \cdot M^p 
  \le 
  \|u_0\|_{L^r_{ul,\rho}(\mathbb{R}^N)}+\frac{T^{1-\alpha p}}{1-\alpha p} M^p.
 \end{split}
\end{equation}
Again by Lemma~\ref{Lemma:2.1} we obtain 
\begin{align*} 
  & 
  t^{\alpha}\|\Phi(u)\|_{L^{pr}_{ul,\rho}(\mathbb{R}^N)} \\ 
  &\ 
  \le 
  C(t^\alpha \rho^{-2\alpha}+1) 
  \|u_0\|_{L^r_{ul,\rho}(\mathbb{R}^N)}+Ct^{\alpha}\int_0^t(\rho^{-2\alpha}+(t-s)^{-\alpha}) \|u(s)\|_{L^{pr}_{ul,\rho}(\mathbb{R}^N)}^p\, ds
  \\
  &\ 
  \le 
  C(t^\alpha \rho^{-2\alpha}+1)
  \|u_0\|_{L^r_{ul,\rho}(\mathbb{R}^N)}+Ct^{\alpha}\int_0^t (\rho^{-2\alpha}+(t-s)^{-\alpha}) s^{-\alpha p}\, ds\cdot M^p
  \\
  &\ 
  \le 
  C(T^\alpha \rho^{-2\alpha}+1)
  \|u_0\|_{L^r_{ul,\rho}(\mathbb{R}^N)}+C\left(\frac{T^{1-\alpha p+\alpha}}{1-\alpha p}\rho^{-2\alpha}+B(1-\alpha,1-\alpha p)T^{1-\alpha p}\right) M^p,
\end{align*}  
where $B(\cdot,\cdot)$ is the beta function and $C>0$ is a constant depending only on $N$, $p$ and $r$. 
Set $M=(2C+4) \|u_0\|_{L^r_{ul,\rho}(\mathbb{R}^N)}$. 
Let $T_0$ be the constant satisfying 
\begin{equation} 
\label{eq:A.2} 
\begin{split} 
& 
C T_0^\alpha \rho^{-2\alpha}  
+(2C+4)^p \|u_0\|_{L^r_{ul,\rho}(\mathbb{R}^N)}^{p-1}\cdot \frac{T_0^{1-\alpha p}}{1-\alpha p} \\ 
& 
\quad +p\, C(2C+4)^p \|u_0\|_{L^r_{ul,\rho}(\mathbb{R}^N)}^{p-1}  
\left(\frac{T_0^{1-\alpha p+\alpha}}{1-\alpha p}\rho^{-2\alpha}+B(1-\alpha,1-\alpha p)T_0^{1-\alpha p}\right)  
=1.  
\end{split} 
\end{equation} 
Then $\Phi$ is a map from $X_{M,T_0}$ to itself. 
Similarly, for $u,v\in X_{M,T}$, we have 
\begin{align*}   
  &
  t^{\alpha}
  \|\Phi(u)-\Phi(v)\|_{L^{pr}_{ul,\rho}(\mathbb{R}^N)}
  \\
  &\quad 
  \le
   Cpt^{\alpha}\int_0^t\left(\rho^{-2\alpha}+(t-s)^{-\alpha}\right)\|u(s)-v(s)\|_{L^{pr}_{ul,\rho}}
   \left(\|u(s)\|_{L^{pr}_{ul,\rho}}^{p-1}+\|v(s)\|_{L^{pr}_{ul,\rho}}^{p-1}\right)   ds
\\
   &\quad 
   \le
   2Cp\left(\frac{T^{1-\alpha p+\alpha}}{1-\alpha p}\rho^{-2\alpha}+B(1-\alpha,1-\alpha p)T^{1-\alpha p}\right) M^{p-1}
   \sup_{0<t<T}t^{\alpha}\|u-v\|_{L^{pr}_{ul,\rho}(\mathbb{R}^N)} 
\end{align*}
for the same constant $C$ as in \eqref{eq:A.2}.
Hence, by \eqref{eq:A.2} wee see that $\Phi$ is a contraction map from $X_{M,T_0}$ to itself. 
Therefore, by the contraction mapping theorem we find a fixed point $u\in X_{M,T_0}$. 
Let $T$ be the maximal existence time such that the fixed point can be found in $X_{M, T}$. 
Then we clearly have $T\ge T_0$, 
and obtain \eqref{eq:2.11} by \eqref{eq:A.2}. 

Next we prove $u\in C((0,T); L^r_{ul,\rho}(\mathbb{R}^N))$. 
Recall that $|u|^p\in L^1(0,T;L^r_{ul,\rho}(\mathbb{R}^N))$ since $u\in X_{M,T}$,  
so $\|u\|_{L^{pr}_{ul,\rho}(\mathbb{R}^N)}^p\le Mt^{-\alpha p}$ for all $t\in (0,T)$. 
This together with the fact $-\alpha p>-1$ implies that 
\begin{equation} 
\label{eq:A.3} 
	u(t)-e^{t\Delta}u_0 =\int_0^t e^{(t-s)\Delta}|u(s)|^{p-1}u(s)\, ds \in C([0,T);\mathcal{L}^r_{ul,\rho}(\mathbb{R}^N)). 
\end{equation} 
On the other hand, since $e^{t\Delta}u_0\in \mathcal{L}^r_{ul,\rho}(\mathbb{R}^N)$, 
we have $e^{t\Delta}u_0\in C((0,T);\mathcal{L}^r_{ul,\rho}(\mathbb{R}^N))$. 
Therefore we obtain $u\in C((0,T);\mathcal{L}^r_{ul,\rho}(\mathbb{R}^N))$. 
We remark that, if $u_0\in \mathcal{L}^r_{ul,\rho}(\mathbb{R}^N)$, then $e^{t\Delta}u_0\in C([0,T);\mathcal{L}^r_{ul,\rho}(\mathbb{R}^N))$, 
so $u\in C([0,T);\mathcal{L}^r_{ul,\rho}(\mathbb{R}^N))$. 

Finally, applying the same iteration argument as in the proof of Proposition~\ref{Proposition:2.3} with $q_1=pr$, instead of $q_1=q$, 
we obtain $L^{\infty}_{loc}((0,T); L^{\infty}(\mathbb{R}^N))$. 
This can be shown similarly, so we give its proof only for the critical case. 
See the argument below. 
Then the standard regularity argument implies that the fixed point $u$ is a classical solution of \eqref{eq:A.1}. 
This solution satisfies
$\|u(t)-e^{t\Delta}u_0\|_{L^r_{ul,\rho}(\mathbb{R}^N)} \to 0$ as $t\to 0$, 
which follows from \eqref{eq:A.3}. 
Thus we complete the proof of Proposition~\ref{Proposition:2.2}. 
\hfill 
\qed 

\vspace{5pt} 

\noindent 
{\bf Proof of Proposition~\ref{Proposition:2.3}.} 
Assume that $r=\frac{N}{2}(p-1)>1$ and $u_0\in \mathcal{L}^r_{ul,\rho}(\mathbb{R}^N)$.
Let $q$ satisfy $\max\{p,r\}<q<pr$ and put $\sigma:=\frac{N}{2}\left(\frac{1}{r}-\frac{1}{q}\right)$.
For any $M>0$ and $T>0$, define 
$$
  Y_{M,T} 
  := \left\{
\begin{array}{l} 
u\in L^{\infty}(0,T; L^{r}_{ul,\rho}(\mathbb{R}^N))\cap L^{\infty}_{loc}((0,T); L^{q}_{ul,\rho}(\mathbb{R}^N)) :\, \vspace{5pt} \\
\displaystyle\sup_{0<t<T}\|u(t)\|_{L^r_{ul,\rho}(\mathbb{R}^N)}\le \|u_0\|_{L^r_{ul,\rho}(\mathbb{R}^N)}+M, \ 
\sup_{0<t<T}t^{\sigma}\|u(t)\|_{L^q_{ul,\rho}(\mathbb{R}^N)}\le M
\end{array}
  \right\},
$$ 
equipped with the metric $d_Y(u,v):=\sup\{t^{\sigma}\|u(t)-v(t)\|_{L^q_{ul,\rho}(\mathbb{R}^N)}:\, 0<t<T\}$.
Then $(Y_{M,T}, d_Y)$ is a complete metric space. 
We show that $\Phi$ is a contraction map from $Y_{M,T}$ to itself for suitable $M>0$ and $T>0$.
Let $u\in Y_{M,T}$. 
Lemma~\ref{Lemma:2.1} shows  
\begin{align*} 
&
\|\Phi(u)\|_{L^r_{ul,\rho}(\mathbb{R}^N)} 
\le 
\|u_0\|_{L^r_{ul,\rho}(\mathbb{R}^N)} 
+
\int_0^t  \|e^{(t-s)\Delta}|u(s)|^p\|_{L^r_{ul,\rho}(\mathbb{R}^N)}\, ds
\\
&\quad 
\le
\|u_0\|_{L^r_{ul,\rho}(\mathbb{R}^N)} 
+
C_1 \int_0^t\left(\rho^{-N\left(\frac{1}{r}-\frac{p}{q}\right)}+(t-s)^{-\frac{N}{2}\left(\frac{1}{r}-\frac{p}{q}\right)}\right)\|u(s)\|^p_{L^q_{ul,\rho}(\mathbb{R}^N)}\, ds
\\
&\quad 
  \le \|u_0\|_{L^r_{ul,\rho}(\mathbb{R}^N)} 
\\
&
\qquad \quad+
 C_1 \left(
 t^{1-\sigma p}\rho^{-N\left(\frac{1}{r}-\frac{p}{q}\right)}
 +
 t^{1-\sigma p-\frac{N}{2}\left(\frac{1}{r}-\frac{p}{q}\right)}
 \int_0^1(1-s)^{-\frac{N}{2}\left(\frac{1}{r}-\frac{p}{q}\right)}s^{-\sigma p}\, 
 ds
 \right)M^p 
\end{align*} 
for some constant $C_1>0$.  
Then, since $-\sigma p>-1$ and $1-\sigma p-\frac{N}{2}\left(\frac{1}{r}-\frac{p}{q}\right)=0$, we have 
\begin{equation}
\label{eq:A.4}%
\sup_{0<t<T}\|\Phi(u)\|_{L^r_{ul,\rho}(\mathbb{R}^N)} 
\le \|u_0\|_{L^r_{ul,\rho}(\mathbb{R}^N)} +C_2 (T^{1-\sigma p}\rho^{-N\left(\frac{1}{r}-\frac{p}{q}\right)}+1)M^p 
\end{equation} 
for some constant $C_2>0$. 
Similarly, we have 
\[
\begin{split}
 &
 t^{\sigma}\|\Phi(u)\|_{L^q_{ul,\rho}(\mathbb{R}^N)}
 \le 
 t^{\sigma}\|e^{t\Delta}u_0\|_{L^q_{ul,\rho}(\mathbb{R}^N)}
 +
 t^{\sigma}\int_0^t\|e^{(t-s)\Delta}|u(s)|^p\|_{L^q_{ul,\rho}(\mathbb{R}^N)}\, ds
 \\
 &\quad
 \lesssim  
 t^{\sigma}\|e^{t\Delta}u_0\|_{L^q_{ul,\rho}(\mathbb{R}^N)}
 +
 t^{\sigma}\int_0^t\left(\rho^{-N\left(\frac{p}{q}-\frac{1}{q}\right)}+(t-s)^{-\frac{N}{2}\left(\frac{p}{q}-\frac{1}{q}\right)}\right)\|u(s)\|^p_{L^q_{ul,\rho}(\mathbb{R}^N)}\, ds
 \\
 &\quad
 \lesssim 
 t^{\sigma}\|e^{t\Delta}u_0\|_{L^q_{ul,\rho}(\mathbb{R}^N)}
\\
&\qquad  +
 \left(
 t^{\sigma+1-\sigma p}\rho^{-N\left(\frac{p}{q}-\frac{1}{q}\right)}
 +
 t^{\sigma+1-\sigma p-\frac{N}{2}\left(\frac{p}{q}-\frac{1}{q}\right)}
 \int_0^1(1-s)^{-\frac{N}{2}\left(\frac{p}{q}-\frac{1}{q}\right)}s^{-\sigma p}\, 
 ds
 \right)M^p.
\end{split}
\]
Then, since 
$-\sigma p>-1$, $-\frac{N}{2}\left(\frac{p}{q}-\frac{1}{q}\right)>-1$ and $\sigma+1-\sigma p-\frac{N}{2}\left(\frac{p}{q}-\frac{1}{q}\right)=0$,  
we obtain 
\begin{equation}
\label{eq:A.5}%
\sup_{0<t<T} t^{\sigma}\|\Phi(u)\|_{L^q_{ul,\rho}(\mathbb{R}^N)}
\lesssim  
\sup_{0<t<T}t^{\sigma}\|e^{t\Delta}u_0\|_{L^q_{ul,\rho}(\mathbb{R}^N)}
+(T^{1+\sigma-\sigma p}\rho^{-N\left(\frac{p}{q}-\frac{1}{q}\right)}+1)M^p.
\end{equation}
It remains to estimate the linear term 
$t^{\sigma}\|e^{t\Delta}u_0\|_{L^q_{ul,\rho}(\mathbb{R}^N)}$.
Since $u_0\in \mathcal{L}^r_{ul,\rho}(\mathbb{R}^N)$, 
there exists a sequence $\{u_n\}\subset BUC(\mathbb{R}^N)$ such that
$\|u_0-u_n\|_{L^r_{ul,\rho}(\mathbb{R}^N)}\to 0$ as $n\to \infty$, which yields 
\begin{align} 
\notag 
 \sup_{0<t<T}t^{\sigma}\|e^{t\Delta}u_0\|_{L^q_{ul,\rho}(\mathbb{R}^N)}
&
\le 
\sup_{0<t<T}t^{\sigma}\|e^{t\Delta}(u_0-u_n)\|_{L^q_{ul,\rho}(\mathbb{R}^N)}
+
\sup_{0<t<T}t^{\sigma}\|e^{t\Delta}u_n\|_{L^q_{ul,\rho}(\mathbb{R}^N)}
\\
\notag 
& \lesssim (T^{\sigma}\rho^{-2\sigma}+1)\|u_0-u_n\|_{L^r_{ul,\rho}(\mathbb{R}^N)}+T^{\sigma}\|u_n\|_{L^q_{ul,\rho}(\mathbb{R}^N)} 
\end{align} 
for sufficiently large $n$ and small $T>0$.
Combining this inequality with \eqref{eq:A.4} and \eqref{eq:A.5}, we obtain 
\[
 \sup_{0<t<T} \|\Phi(u)\|_{L^r_{ul,\rho}}\le \|u_0\|_{L^r_{ul,\rho}}+M
\quad \text{and}\quad
 \sup_{0<t<T} t^{\sigma}\|\Phi(u)\|_{L^q_{ul,\rho}}
\le M
\]
for suitable $M>0$ and $T>0$.
Similarly, we get 
\begin{equation}
\label{eq:A.6}
  \begin{split} 
    &
    \sup_{0<t<T} t^{\sigma}\|\Phi(u)-\Phi(v)\|_{L^q_{ul,\rho}(\mathbb{R}^N)} 
	\\ 
	&
	\quad \qquad 
	\lesssim 
	\left( T^{1+\sigma-\sigma p}\rho^{-N\left(\frac{p}{q}-\frac{1}{q}\right)}+1 \right) M^{p-1}
    \sup_{0<t<T} t^{\sigma}\|u-v\|_{L^q_{ul,\rho}(\mathbb{R}^N)} 
  \end{split} 
\end{equation}
for all $u,v\in Y_{M,T}$.
This proves that $\Phi$ is a contraction map from $Y_{M,T}$ to itself
for sufficiently small $M>0$ and $T>0$. 
Thus, by the contraction mapping theorem we can find a fixed point $u\in Y_{M,T}$. 

Let us prove that  $u\in Y_{M,T}$ satisfies
$u\in C((0,T);\mathcal{L}^q_{ul,\rho}(\mathbb{R}^N))$ and
$t^{\sigma}\|u\|_{L^q_{ul,\rho}(\mathbb{R}^N)}\to 0$ as $t\to 0$.
Put 
\[
K:=
Y_{M,T} 
\cap 
\left\{
u\in C((0,T);\mathcal{L}^q_{ul,\rho}(\mathbb{R}^N))
:\, \lim_{t\to 0} t^{\sigma}\|u\|_{L^q_{ul,\rho}(\mathbb{R}^N)}=0
 \right\}.
\]
Then $K$ equipped with a metric $d_Y$ is a complete metric space.
We now prove that $\Phi$ is a map from $K$ to $K$
(then the contraction mapping argument works in $K$, therefore the fixed point $u$ belongs to $K$). 
Assume that $u_0\in \mathcal{L}^r_{ul,\rho}(\mathbb{R}^N)$.
By the smoothing effect of the heat semigroup~$e^{t\Delta}$ 
we have $e^{t\Delta}u_0\in \mathcal{L}^q_{ul,\rho}(\mathbb{R}^N)\subset \mathcal{L}^r_{ul,\rho}(\mathbb{R}^N)$ for $t>0$.
This and Lemma~\ref{Lemma:2.2} imply that $\|e^{(t+h)\Delta}u_0-e^{t\Delta}u_0\|_{L^q_{ul,\rho}(\mathbb{R}^N)}\to 0 $ as $h\to 0$ for all $t>0$.
Moreover, by the same argument as in the above argument 
we have $t^{\sigma}\|e^{t\Delta}u_0\|_{L^q_{ul,\rho}(\mathbb{R}^N)}\to 0$ as $t\to 0$. 
Thus we obtain $e^{t\Delta}u_0\in K$. 
Since $K\cap C((0,T); BUC(\mathbb{R}^N))$ is dense in $K$ equipped with $d_Y$,
there exists a sequence $\{u_n\}\subset K\cap C((0,T); BUC(\mathbb{R}^N))$ such that $d_Y(u,u_n)\to 0$ as $n\to \infty$.
This together with the facts $u,u_n \in  Y_{M,T}$ and \eqref{eq:A.6} yields $d_Y(\Phi(u),\Phi(u_n))\to 0$ as $n\to \infty$. 
Furthermore, we have $\Phi(u_n)\in K$, which follows from the fact $e^{t\Delta}u_0\in K$ and $u_n\in C((0,T); BUC(\mathbb{R}^N))$. 
Therefore, since $K$ is a complete metric space equipped with $d_Y$, we obtain $\Phi(u)\in K$. 

We are in position to prove $u\in C([0,T);\mathcal{L}^r_{ul,\rho}(\mathbb{R}^N))$. 
It suffices to prove 
\begin{equation} 
\label{eq:A.7}%
\lim_{t\to 0}\|u(t)-u_0\|_{L^r_{ul,\rho}(\mathbb{R}^N)}=0, 
\end{equation} 
since 
$K\subset C((0,T);\mathcal{L}^r_{ul,\rho}(\mathbb{R}^N))$.
We clearly see 
that $\|e^{t\Delta}u_0-u_0\|_{L^r_{ul,\rho}(\mathbb{R}^N)}\to 0$ as $t\to 0$ by $u_0\in \mathcal{L}^r_{ul,\rho}(\mathbb{R}^N)$ and Lemma~\ref{Lemma:2.2}. 
Moreover, we have 
\begin{align*} 
  & 
  \|u(t)-e^{t\Delta}u_0\|_{L^r_{ul,\rho}(\mathbb{R}^N)} \\ 
  &\quad 
  \le 
  \int_0^t \|e^{(t-s)\Delta}|u(s)|^p\|_{L^r_{ul,\rho}(\mathbb{R}^N)}\, ds 
  \\
  &\quad 
  \lesssim  
  \int_0^t
 \left(
  \rho^{-N\left(\frac{p}{q}-\frac{1}{r}\right)}+(t-s)^{-\frac{N}{2}\left(\frac{p}{q}-\frac{1}{r}\right)}
 \right)s^{-\frac{N}{2}\left(\frac{1}{r}-\frac{1}{q}\right)p}\, ds
 \cdot \sup_{0<s<t}s^{\sigma p}\|u(s)\|_{L^q_{ul,\rho}(\mathbb{R}^N)}^p 
 \\
 &\quad 
 \lesssim 
 \left(
  t^{1-\frac{N}{2}\left(\frac{1}{r}-\frac{1}{q}\right)p}\rho^{-N\left(\frac{p}{q}-\frac{1}{r}\right)}+1
 \right)
  \sup_{0<s<t}s^{\sigma p}\|u(s)\|_{L^q_{ul,\rho}(\mathbb{R}^N)}^p
\to 0\ \ \text{as}\ \ t\to 0,
\end{align*} 
since $u\in K$ and  
$$ 
-\frac{N}{2}\left(\frac{1}{r}-\frac{1}{q}\right)p>-1, \qquad 
-\frac{N}{2}\left(\frac{p}{q}-\frac{1}{r}\right)-\frac{N}{2}\left(\frac{1}{r}-\frac{1}{q}\right)p=-1. 
$$
Thus we have \eqref{eq:A.7}, and so $u\in C([0,T);\mathcal{L}^r_{ul,\rho}(\mathbb{R}^N))$. 

Finally we prove that the fixed point $u\in C([0,T); L^r_{ul,\rho}(\mathbb{R}^N))$ is smooth so is a classical solution of \eqref{eq:A.1}. 
To this end, it suffices to prove $u\in L^{\infty}_{loc}((0,T); L^{\infty}(\mathbb{R}^N))$. 
Let $n\in \mathbb{N}$ and choose $\{q_n\}$ such that 
\begin{equation}
\notag 
q_1:=q,\ \ q_{n+1}\ge q_n\ \ \text{with}\ \ \frac{N}{2}\left(\frac{p}{q_{n}}-\frac{1}{q_{n+1}}\right)<1.
\end{equation}
A simple calculation shows that
\[
 \frac{1}{q_n}>p^{n-1}\left(\frac1q-\frac{p}{r}\right)+\frac{p}{r}.
\]
Since 
$\frac1q-\frac{p}{r}<\frac1q-\frac{1}{r}<0$,
there exists $n_0\in \mathbb{N}$
such that
$\frac{Np}{2q_{n_0-1}}<1$.
Then we redefine $q_{n_0}:=\infty$ and consider the sequence $\{q_n\}_{n=1}^{n_0}$.
Now fix $\epsilon>0$.
Since $u\in Y_{M,T}$, 
we have $u\in L^{\infty} ( [{\epsilon}/{n_0},T] ; L^{q_1}_{ul,\rho}(\mathbb{R}^N) )$ and
$$
 u(t+\epsilon/{n_0})=e^{t\Delta}\left(u(\epsilon/{n_0})\right)+\int_0^te^{(t-s)\Delta}|u(s+\epsilon/{n_0})|^{p-1}u(s+\epsilon/{n_0})\, ds, 
$$
which implies that there exists a constant $C(\epsilon)>0$ such that
\begin{align*} 
& \left\|u\left(t+\frac{\epsilon}{n_0}\right)\right\|_{L^{q_2}_{ul,\rho}(\mathbb{R}^N)} \\ 
& \quad \lesssim 
 \left(\rho^{-N\left(\frac{1}{q_2}-\frac{1}{q_1}\right)} 
 +t^{-\frac{N}{2}\left(\frac{1}{q_2}-\frac{1}{q_1}\right)}\right)\left\|u\left(\frac{\epsilon}{n_0}\right)\right\|_{L^{q_1}_{ul,\rho}(\mathbb{R}^N)} \\ 
& \quad \qquad +
 \int_0^t\left(\rho^{-N\left(\frac{p}{q_1}-\frac{1}{q_2}\right)} 
 +(t-s)^{-\frac{N}{2}\left(\frac{p}{q_1}-\frac{1}{q_2}\right)}\right)\left\|u\left(\cdot+\frac{\epsilon}{n_0}\right)\right\|_{L_{ul,\rho}^{q_1}(\mathbb{R}^N)}^p ds 
 \le C(\epsilon)
\end{align*} 
for all $t\in[\epsilon/{n_0},T-\epsilon/{n_0}]$.
This proves $u\in L^{\infty}([2\epsilon/{n_0},T]; L^{q_2}(\mathbb{R}^N))$.
Iterating above argument $n_0$ times,
we obtain $u\in L^{\infty}([\epsilon,T]; L^{q_{n_0}}(\mathbb{R}^N))$.
Since $q_{n_0}=\infty$ and $\epsilon>0$ is arbitrary, we have $u\in L^{\infty}_{loc}((0,T); L^{\infty}(\mathbb{R}^N))$. 
Thus we complete the proof of Proposition~\ref{Proposition:2.3}. 
\hfill 
\qed

\medskip 

\noindent 
{\bf Acknowledgments.} 
This work was partially funded by JSPS KAKENHI (grant number 15K17573 and 15K17575)
and Sumitomo Fundation (grant number 140825).


\begin{thebibliography}{xx} 
\bibitem{AD} 
D.~Andreucci and E.~DiBenedetto,
{\it On the Cauchy problem and initial traces for a class of evolution equations with strongly nonlinear sources}, 
Ann. Scuola Norm. Sup. Pisa Cl. Sci. {\bf 18} (1991), 363--441. 

\bibitem{ABCD} 
J.~M.~Arriera, A.~Rodriguez-Bernal, J.~W.~Cholewa and T.~Dlotko, 
{\it Linear parabolic equations in locally uniform spaces}, 
Math. Models  Methods Appl. Sci. {\bf 14} (2004), 253--293. 

\bibitem{BP} 
P.~Baras and M.~Pierre, 
{\it Crit\`{e}re d'existence de solutions positives pour des \'{e}quations semi-lin\'{e}aires non monotones}, 
Ann. Inst. H. Poincar\'{e} Anal. Non Lin\'{e}aire {\bf 2} (1985), 185--212. 

\bibitem{BC}
H.~Brezis and T.~Cazenave,
{\it A nonlinear heat equation with singular initial data},
J. Anal. Math. {\bf 68} (1996), 277--304. 

\bibitem{F} 
Y.~Fujishima, 
{\it Blow-up set for a superlinear heat equation and pointedness of the initial data}, 
Discrete Continuous Dynamical Systems A {\bf 34} (2014), 4617--4645. 

\bibitem{G} 
Y.~Giga, 
{\it Solutions for semilinear parabolic equations in $L^p$ and regularity of weak solutions of the Navier-Stokes system}, 
J. Differential Equations {\bf 62} (1986), 415--421. 

\bibitem{GV} 
J.~Ginibre and G.~Velo, 
{\it The Cauchy problem in local spaces for the complex Ginzburg-Landau equation II. Contraction methods}, 
Commun. Math. Phys. {\bf 187} (1997), 45--79. 

\bibitem{HW} 
A.~Haraux and F.~B.~Weissler, 
{\it Nonuniqueness for a semilinear initial value problem}, 
Indiana Univ. Math. J. {\bf 31} (1982), 167--189. 

\bibitem{IJMS} 
S.~Ibrahim, R.~Jrad, M.~Majdoub and T.~Saanouni, 
{\it Local well posedness of a 2D semilinear heat equation}, 
Bull. Belg. Math. Soc. Simon Stevin {\bf 21} (2014), 535--551.  

\bibitem{I} 
N.~Ioku, 
{\it The Cauchy problem for heat equations with exponential nonlinearity}, 
J. Differential Equations {\bf 251} (2011), 1172--1194. 

\bibitem{IRT} 
N.~Ioku, B.~Ruf and E.~Terraneo, 
{\it Existence, non-existence, and uniqueness for a heat equation with exponential nonlinearity in $\mathbb{R}^N$}, 
Math. Phys. Anal. Geom. {\bf 18} (2015), 18:29.

\bibitem{IKS} 
K.~Ishige, T.~Kawakami and M.~Sier\.{z}\c{e}ga, 
{\it Supersolutions of parabolic systems with power nonlinearities}, 
J. Differential Equations {\bf 260} (2016), 6084--6107. 

\bibitem{IS} 
K.~Ishige and R.~Sato, 
{\it Heat equation with a nonlinear boundary condition and uniformly local $L^r$ spaces}, 
Discrete Continuous Dynamical Systems A {\bf 36}, 2627--2652. 

\bibitem{LSU}
O.~A.~Ladyzenskaja, N.~S.~Solonikov and N.~N.~Ural'ceva, 
``Linear and quasilinear equations of parabolic type'', 
Amer. Math. Soc., Providence, 1968. 

\bibitem{LRSV} 
R.~Laister, J.C.~Robinson, M.~Sier\.{z}\c{e}ga and A.~Vidal-L\'{o}pez, 
A complete characterization of local existence of semilinear heat equations in Lebesgue spaces, 
Ann. Inst. H. Poincar\'{e} Anal. Non Lin\'{e}aire, in press. 

\bibitem{MT} 
Y. Maekawa and Y. Terasawa,
{\it The Navier-Stokes equations with initial data in uniformly local $L^p$ spaces}, 
Differential Integral Equations {\bf 19} (2006), 369--400. 

\bibitem{NO} 
M.~Nakamura and T.~Ozawa, 
{\it Nonlinear Schr\"{o}dinger equations in the Sobolev space of critical order}, 
J. Funct. Anal. {\bf 155} (1998), 364--380. 

\bibitem{NS} 
W.-M.~Ni and P.~Sacks, 
{\it Singular behavior in nonlinear parabolic equations}, 
Trans. Amer. Math. Soc. {\bf 287} (1985), 657--671. 

\bibitem{QS}
P. Quittner and P. Souplet, 
``Superlinear Parabolic Problems, Blow-up, Global Existence and Steady States'', 
Birkh\"auser Advanced Texts: Basler Lehrb\"ucher
Birkh\"auser Verlag, Basel, 2007. 

\bibitem{RS} 
J.~C.~Robinson and M.~Sier\.{z}\c{e}ga, 
{\it Supersolutions for a class of semilinear heat equations}, 
Rev. Mat. Complut. {\bf 26} (2013), 341--360.  

\bibitem{RT} 
B.~Ruf and E.~Terraneo, 
{\it The Cauchy problem for a semilinear heat equation with singular initial data}, 
Progr. Nonlinear Differential Equations Appl. {\bf 50} (2002), 295--309. 

\bibitem{S} 
D.~H.~Sattinger, 
{\it Monotone methods in nonlinear elliptic and parabolic boundary value problems}, 
Indiana Univ. Math. J. {\bf 21} (1972), 979--1000. 

\bibitem{STW} 
S.~Snoussi, S.~Tayachi and F.~B.~Weissler, 
{\it Asymptotically self-similar global solutions of a semilinear parabolic equation with a nonlinear gradient term}, 
Proc. Roy. Soc. Edinburgh Sect. A {\bf 129} (1999), 1291--1307. 

\bibitem{Ta} 
S.~Tayachi, 
{\it Forward self-similar solutions of a semilinear parabolic equation with a nonlinear gradient term}, 
Differential Integral Equations {\bf 9} (1996), 1107--1117. 

\bibitem{Te} 
E.~Terraneo, 
{\it Non-uniqueness for a critical non-linear heat equation}, 
Comm. Partial Differential Equations {\bf 27} (2002), 185--218. 

\bibitem{W1}
F.~B.~Weissler, 
{\it Local existence and nonexistence for semilinear parabolic equations in $L^p$}, 
Indiana Univ. Math. J. {\bf 29} (1980), 79--102.


\end{thebibliography}
\end{document}